\documentclass[12pt]{article}

\usepackage{subcaption}

\usepackage{float}

\usepackage{multirow}

\usepackage{tikz}

\usepackage{algorithm}
\usepackage{algpseudocode}
\algrenewcommand\algorithmicrequire{\textbf{Input:}}
\algrenewcommand\algorithmicensure{\textbf{Output:}}

\usetikzlibrary{calc,backgrounds,arrows,matrix}

\usepackage{enumerate}
\usepackage{listings}

\usepackage{color}
\definecolor{lightblue}{rgb}{0,0.2,0.5}

\usepackage[colorlinks=true, urlcolor=lightblue,linkcolor=lightblue, citecolor=lightblue]{hyperref}

\usepackage{amssymb,amsmath}

\usepackage{color}
\usepackage{mathrsfs}

\usepackage{mathtools}
\usepackage{accents}

\DeclareMathAlphabet{\eufrak}{U}{}{}{}
\DeclarePairedDelimiter\floor{\lfloor}{\rfloor}
\SetMathAlphabet\eufrak{normal}{U}{euf}{m}{n}
\SetMathAlphabet\eufrak{bold}{U}{euf}{b}{n}

\DeclareMathOperator*{\argmin}{arg\,min}

\allowdisplaybreaks

 \def\qu{{\mathord{\mathbb Z}}}

 \def\inte{{\mathord{\mathbb R}}}

 \def\inte{{\mathord{\mathbb N}}}

 \def\sZZ{{\rm Z\kern-.45em{}Z}}

 \def\sQQ{{\kern 0.27em \vrule height1.45ex width0.03em depth0em
           \kern-0.30em \rm Q}}
 \def\qu{{\mathchoice
         {\sQQ}
         {\sQQ}
   {\kern 0.225em \vrule height1.05ex width0.025em depth0em \kern-0.25em \rm Q}
   {\kern 0.180em \vrule height0.78ex width0.020em depth0em \kern-0.20em \rm Q}
         }}
 \def\sGG{{\kern 0.27em \vrule height1.45ex width0.03em depth0em
           \kern-0.30em \rm G}}
 \def\gg{{\mathchoice
         {\sGG}
         {\sGG}
   {\kern 0.225em \vrule height1.05ex width0.025em depth0em \kern-0.25em \rm G}
   {\kern 0.180em \vrule height0.78ex width0.020em depth0em \kern-0.20em \rm G}
         }}

 \newtheorem{prop}{Proposition}[section]
 
 \newtheorem{definition}[prop]{Definition}

 \newtheorem{remark}[prop]{Remark}

\numberwithin{equation}{section}

\newcommand{\abs}[1]{\lvert#1\rvert}

\newcommand{\re}{\mathrm{e}}

 \newcounter{hyp}
\newenvironment{Proof}{\removelastskip\par\medskip \noindent{\em Proof.} \rm}{\penalty-20\null\hfill$\square$\par\medbreak}

\def\bprf{\begin{Proof}}
\def\nprf{\end{Proof}}
\def\bdes{\begin{description}}
\def\ndes{\end{description}}

\newtheorem{thm}{Theorem}[section]

\def\bdef{\begin{defn}}
\def\ndef{\end{defn}}
\def\bthm{\begin{thm}}
\def\nthm{\end{thm}}
\def\bprop{\begin{prop}}
\def\nprop{\end{prop}}
\def\brmk{\begin{remark}}
\def\nrmk{\end{remark}}
\def\bexa{\begin{exa}}
\def\nexa{\end{exa}}
\def\blem{\begin{lem}}
\def\nlem{\end{lem}}
\def\bcor{\begin{cor}}
\def\ncor{\end{cor}}
\def\bexe{\begin{exe}}
\def\nexe{\end{exe}}

\newcommand{\E}{\mathbb{E}}

\newcommand{\nn}{\mathbb{N}}
\newcommand{\real}{\mathbb{R}}

\newcommand{\domain}{D}

\def\og{\leavevmode\raise.3ex
     \hbox{$\scriptscriptstyle\langle\!\langle$~}}
\def\fg{\leavevmode\raise.3ex
     \hbox{~$\!\scriptscriptstyle\,\rangle\!\rangle$}~}

\title{\Huge
 A deep branching solver for 
 fully nonlinear partial differential equations
}

\author{
 Jiang Yu Nguwi\footnote{\href{mailto:nguw0003@e.ntu.edu.sg}{nguw0003@e.ntu.edu.sg}
 }
 \qquad
 Guillaume Penent\footnote{\href{mailto:PENE0001@e.ntu.edu.sg}{pene0001@e.ntu.edu.sg}}
 \qquad Nicolas Privault\footnote{
\href{mailto:nprivault@ntu.edu.sg}{nprivault@ntu.edu.sg}
 }
 \\
  \small
Division of Mathematical Sciences
\\
\small
School of Physical and Mathematical Sciences
\\
\small
Nanyang Technological University
\\
\small
21 Nanyang Link, Singapore 637371
}

\allowdisplaybreaks

\usepackage{empheq}
\usepackage{bm}

\usepackage{bbm}

\makeatletter
\newcommand*\rel@kern[1]{\kern#1\dimexpr\macc@kerna}
\newcommand*\widebar[1]{
  \begingroup
  \def\mathaccent##1##2{
    \rel@kern{0.8}
    \overline{\rel@kern{-0.8}\macc@nucleus\rel@kern{0.2}}
    \rel@kern{-0.2}
  }
  \macc@depth\@ne
  \let\math@bgroup\@empty \let\math@egroup\macc@set@skewchar
  \mathsurround\z@ \frozen@everymath{\mathgroup\macc@group\relax}
  \macc@set@skewchar\relax
  \let\mathaccentV\macc@nested@a
  \macc@nested@a\relax111{#1}
  \endgroup
}
\makeatother

\begin{document}
\maketitle

\baselineskip0.6cm

\vspace{-0.6cm}

\begin{abstract}
We present a multidimensional deep learning implementation of a stochastic branching algorithm for the numerical solution of fully nonlinear PDEs. This approach is designed to tackle functional nonlinearities involving gradient terms of any orders, by combining the use of neural networks with a Monte Carlo branching algorithm. In comparison with other deep learning PDE solvers, it also allows us to check the consistency of the learned neural network function. Numerical experiments presented show that this algorithm can outperform deep learning approaches based on backward stochastic differential equations or the Galerkin method, and provide solution estimates that are not obtained by those methods in fully nonlinear examples.
\end{abstract}



\baselineskip0.7cm

\parskip-0.1cm

\section{Introduction}
This paper is concerned with the numerical solution of
 fully nonlinear partial differential equations (PDEs) of the form
\begin{equation}
\label{eq:main pde}
\begin{cases}
  \displaystyle
  \partial_t u(t,x) + \frac{1}{2}\Delta u(t,x)
        + f\big(\partial_{\lambda^1}u(t,x) , \ldots , \partial_{\lambda^n}u(t,x)\big) = 0,
  \medskip
  \\
u(T,x) = \phi (x), \qquad (t,x) = (t,x_1, \ldots, x_d) \in [0,T] \times \real^d,
\end{cases}
\end{equation}
$d\geq 1$, where
$    \Delta  = \sum\limits_{i = 1}^d \partial^2 / \partial x^2_i$
            is the standard $d$-dimensional Laplacian,
$\partial_t u (t,x) = \partial u(t,x) / \partial t$,
and $f$ is a smooth function of the derivatives
$$
\partial_{\lambda^i} u(t,x)
=
\frac{\partial^{\lambda_1^i}}{\partial x_1}
\cdots
\frac{\partial^{\lambda_d^i}}{\partial x_d}
u(t,x_1,\ldots , x_d), \qquad (x_1,\ldots , x_d)\in \real^d,
$$
$\lambda^i = \left(\lambda^i_1, \dots, \lambda^i_d\right) \in \nn^d$,
$i = 1, \ldots, n$.
As is well known, standard numerical schemes for solving \eqref{eq:main pde}
by e.g. finite differences or finite elements
suffer from the curse of dimensionality as
 their computational cost grows exponentially with the dimension $d$.

\medskip

The deep Galerkin method (DGM) has been developed in \cite{sirignano2018dgm}
for the numerical solution of \eqref{eq:main pde} by training a neural
network function $v(t, x)$ using the loss function
\begin{equation}
\label{eq:dgm loss function}
\left(\partial_t v(t,x) + \frac{1}{2}\Delta v(t,x)
+
f\big(\partial_{\lambda^1}v(t,x) , \ldots , \partial_{\lambda^n}v(t,x)\big)
\right)^2
+ \left(v(T,x) - \phi (x)\right)^2. 
\end{equation}
{
 See \cite{lyu} for recent improvements of the DGM 
 using deep mixed residuals (MIM) with numerical applications to linear PDEs,
 and \cite{hou} for the blocked residual connection method (DLBR)
 applied to a linear (generalized) Black-Scholes equation.}
\label{p2}
 
\medskip
 
 On the other hand, probabilistic schemes
provide a promising direction to overcome the curse of
dimensionality.
For example, when $f(u(t, x)) = r u(t, x)$ does not involve any derivative of $u$,
the solution of the PDE
$$
\begin{cases}
  \displaystyle
  \partial_t u(t,x) + \frac{1}{2}\Delta u(t,x)
        + r u(t,x) = 0,
  \medskip
  \\
u(T,x) = \phi (x), \qquad (t,x) = (t,x_1, \ldots, x_d) \in [0,T] \times \real^d,
\end{cases}
$$
 admits the probabilistic representation
$$
    u(0, x) = \re^{rT} \E [ \phi (x + W_T) ],
$$
where $(W_t)_{t \ge 0}$ is a standard Brownian motion.
This method can be implemented on a bounded domain $\domain \subset \real^d$
 based on the universal approximation theorem
and the $L^2$ minimality property
$$
    u(0, \cdot) = \inf\limits_v \E \big[
            \big(\re^{rT} \phi(X + W_T) - v(X)\big)^2 \big],
$$
where $X$ is a uniform random vector on $\domain$
and the infimum in $v$ is taken over a neural functional space. 

\medskip

Probabilistic representations for the solutions of
first order nonlinear PDEs 
 can also be obtained
by
 representing $u(t,x)$ as
 $u(t,x) =
 Y_t^{t,x}$, $(t,x)\in [0,T]\times \real$,
 where $(Y_s^{t,x})_{t \leq s \leq T}$ is the
 solution of a backward stochastic differential equation (BSDE),
 see \cite{peng2}, \cite{pardouxpeng}.
 {
 The BSDE method has been implemented in
 \cite{han2018solving} using a deep learning
 algorithm in the case where $f$ depends on the first order derivative,
 i.e.\ $\lambda^i_1 + \cdots + \lambda^i_d \le 1$,
 $1 \le i \le n$, see also \cite{hure2019some}
 for recent improvements. 
 The BSDE method extends to second order fully nonlinear PDEs 
 by the use of second order backward stochastic differential
 equations, see e.g. \cite{touzi}, \cite{soner},
 and \cite{han2018solvingarxiv,beck},
 and \cite{germain}, \cite{lefebvre},
 for deep learning implementations. 
 However, this approach does not apply to nonlinearities
 in gradients of order strictly greater than two, see 
 Examples~e) and f) below. 
}

\medskip 

\label{p3}
{
 Numerical solutions of semilinear PDEs have also
 been obtained by the multilevel Picard method (MLP), see  
 \cite{hutzenthaler-mlp0,hutzenthaler-mlp2,hutzenthaler-mlp3,hutzenthaler-mlp1},
 with numerical experiments provided in \cite{hutzenthaler-mlp4}. 
 However, this approach is currently restricted to first order gradient
 nonlinearities,
 similarly to the deep splitting algorithm of \cite{beck2019deep}. 
 In addition, the main use of the MLP and deep splitting methods 
 is to provide pointwise estimates, whereas this paper
 focuses on functional estimation of solutions using neural networks.}

\medskip
 
{In this context, 
  the use of stochastic branching diffusion mechanisms \cite{skorohodbranching}, \cite{inw}, represents an alternative to the DGM and BSDE methods,
  see \cite{hpmckean}
  for an application to the Kolmogorov-Petrovskii-Piskunov (KPP) equation,
 \cite{chakraborty} for existence of solutions
 of parabolic PDEs with 
 power series nonlinearities, \cite{henry-labordere2012} 
 for more general PDEs with polynomial nonlinearities,
 and \cite{labordere2}
 for an application to semilinear and higher-order hyperbolic PDEs.}
 This approach has been applied in e.g. \cite{lm}, \cite{labordere} to polynomial gradient nonlinearities,
see also \cite{fahim}, \cite{tanxiaolu}, \cite{guowenjie}, 
\cite{huangshuo} for finite difference schemes combined with Monte Carlo estimation for fully nonlinear PDEs with gradients of order up to $2$.
  
\medskip

 Extending such approaches to nonlinearities involving gradients 
 of order greater than two involves technical difficulties
 linked to the integrability of the Malliavin-type weights used in
 repeated integration by parts argument, see page~199 of \cite{labordere}. 
 {Such higher order nonlinearities are also not covered by 
   multilevel Picard \cite{hutzenthaler-mlp4} and
   deep splitting \cite{beck2019deep} methods,
 or by BSDE methods \cite{han2018solving,beck},
 which are limited to first and second order gradients,
 respectively.}
  
\medskip

 In \cite{penent2022fully}, a stochastic branching method 
 that carries information on (functional) nonlinearities along a random tree
 has been introduced, with the aim of providing
 Monte Carlo schemes for the numerical solution of
 fully nonlinear PDEs with gradients of arbitrary orders.

\medskip

 In this paper, we present a deep learning implementation of the
 method of \cite{penent2022fully} using Monte Carlo sampling,
the law of large numbers, and the universal approximation theorem.
 Our approach to the numerical solution of the PDE \eqref{eq:main pde} is based on the following steps:
\begin{enumerate}[i)]
    \item The solution of PDE \eqref{eq:main pde}
            is written as the conditional expectation of a functional of a random coding tree
            via the fully nonlinear Feynman-Kac formula
            Theorem~1 in \cite{penent2022fully},
            see \eqref{eq:feynman kac} below.
    \item The conditional expectation is approximated
            by a neural network function
            through the $L^2$-minimality property
            and the universal approximation theorem.
\end{enumerate}
 We start by testing our method on the Allen-Cahn equation \eqref{eq:example 1},
for which we report a performance comparable to that of the deep BSDE and
deep Galerkin methods, see Figure~\ref{figexpl}. 
This is followed by an example \eqref{eq:example 2}
involving an exponential nonlinearity without gradient
term, in which our method outperforms 
the deep Galerkin method 
and performs comparably to deep BSDE method in dimension $d=5$,
see Figure~\ref{fig:example 2}. 
{We also consider a multidimensional Burgers equation
 \eqref{burgers0-00} 
 for which the deep branching method
 is more stable than the deep Galerkin and deep BSDE methods
 in dimension $d=15$, see Figure~\ref{figexpl-2}}. 
Next, we consider a Merton problem \eqref{eq:example 7} to which
the deep Galerkin method does not apply since
its loss function involves
a division by the second derivative of the neural network function.
{We also note that the deep branching method
  overperforms the deep BSDE method in this case,
see Figure~\ref{fig:example 7}}. 
Finally, we consider higher order functional gradient nonlinearities
in Equations~\eqref{eq:example 6} and \eqref{eq:example 5},
{to which the deep BSDE, multilevel Picard
  and deep splitting methods do not apply}. In those cases, our method
also outperforms the deep Galerkin method in both dimensions $d=1$ and $d=5$,
see Figures~\ref{fig:example 6} and \ref{fig:example 5}. 

\medskip

We also note that since the deep branching method is based on a direct Monte Carlo estimation, it allows for checking the consistency between
the Monte Carlo samples and the learned neural network function,
which is not possible with the deep Galerkin method and deep BSDE methods, 
see Figure~\ref{fig:example 7 debug}.

\medskip

\label{p4}{
 Our algorithm,
 similarly to other branching diffusion methods,
 suffers from a time explosion phenomenon
 due to the use of a branching process.
 Nevertheless, our method can perform better
  than the deep Galerkin and 
  deep BSDE methods in small time
  and in higher dimensions,
  see Figure~\ref{figexpl} 
 for the Allen-Cahn equation and
 Figure~\ref{figexpl-2} for the Burgers equation.}
 
\medskip

\label{p-4}
{Other approaches to the solution of
  evolution equations 
  by carrying information on nonlinearities along trees
  include \cite{butcher1963}, see also 
  Chapters~4-6 of \cite{deuflhard} and 
  \cite{mclachlan}
  for ordinary differential equations (ODEs),
  with applications ranging from geometric numerical integration to 
  stochastic differential equations, see for instance \cite{ehairer}
  and references therein.
 On the other hand, the stochastic branching method does not use
 series truncations and it can be used to estimate an infinite series,
 see \cite{penent4} for an application to ODEs.~ 
}

\medskip 

This paper is organized as follows.
 The extension of the fully nonlinear Feynman-Kac formula of \cite{penent2022fully} to a multidimensional setting is presented in Section~\ref{subsec:feynman kac}, 
 and the deep learning algorithm is described in Section~\ref{sec:main ideas}. 
Section~\ref{sec:numerical examples} presents numerical examples in which
 our method can outperform the deep BSDE and deep Galerkin methods. 

 \medskip 

 The Python codes and numerical experiments run in this paper are available at 
 \\
 \centerline{\url{https://github.com/nguwijy/deep_branching}.}

\subsubsection*{Notation}
\label{subsec:notation}
We denote by $\nn = \{0, 1, 2, \dots\}$ the set of natural numbers, and 
 let $C^{0, \infty}([0, T] \times \real^d)$
be the set of functions $u:[0, T] \times \real^d \to \real$
such that $u(t, x)$ is continuous in the variable $t$ 
and infinitely $x$-differentiable.
For a vector $x = (x_1, \ldots, x_d)^\top \in \real^d$,
we let $\abs{x} = \sum\limits_{i=1}^d \abs{x_i}$, 
and let 
$\bm{1}_p$ be the vector of $1$ at position $p$ and $0$ elsewhere.
We also consider the linear order $\prec$ on $\real^d$ such that
$(k_1, \dots, k_d) = k  \prec l = (l_1, \dots, l_d)$ 
if one of the following holds:
\begin{enumerate}[i)]
    \item $\abs{k} < \abs{l}$;
    \item $\abs{k} = \abs{l}$ and $k_1 < l_1$;
    \item $\abs{k} = \abs{l}$, $k_1 = l_1, \dots k_i = l_i$,
            and $k_{i+1} < l_{i+1}$ for some $1 \le i < d$.
\end{enumerate}

\section{Fully nonlinear Feynman-Kac formula}
\label{subsec:feynman kac}
In this section we extend the construction of
\cite{penent2022fully} to the case of
 multidimensional PDEs of the form
 \begin{equation}
   \label{eq:112}
\begin{cases}
  \displaystyle
  \partial_t u(t,x) + \frac{1}{2}\Delta u(t,x)
  + f
  \big(\partial_{\lambda^1}u(t,x) , \ldots , \partial_{\lambda^n}u(t,x)\big)
  = 0,
  \medskip
  \\
u(T,x) = \phi(x), \qquad (t,x) = (t,x_1, \ldots, x_d) \in [0,T] \times \real^d,
\end{cases}
\end{equation}
 where
 $\lambda^i = \left(\lambda_1^i, \ldots, \lambda_d^i\right) \in \mathbb{N}^d$,
 $i=1,\ldots, n$, with the integral formulation
\begin{equation}
\nonumber
    u (t,x)  =  \int_{\real^d} \varphi (T-t,y-x)\phi(y)dy
   +
   \int_t^T \int_{\real^d} \varphi (s-t,y-x)
   f
   \big(\partial_{\lambda^1}u(t,y) , \ldots , \partial_{\lambda^n}u(t,x)\big)
      dy ds,
\end{equation}
 where $\varphi (t,x) := \re^{- x^2 / ( 2t) } / {\sqrt{2\pi t} }$, and
 $(t,x) \in [0,T]\times \real^d$.
 We refer to e.g. Theorem~1.1 in \cite{krylov1983}
 for sufficient conditions for existence and uniqueness of
 smooth solutions to such fully nonlinear PDEs in the second order case.
 Our fully nonlinear Feynman-Kac formula \cite{penent2022fully} relies on the
 construction of a branching coding tree, based on the definition of
 a set $\mathcal{C}$ of codes and its associated mechanism $\mathcal{M}$.
In what follows, we use the notation
$$(a_1 , \ldots , a_n) \cup (b_1 , \ldots , b_m)
:=
(a_1 , \ldots , a_n, b_1 , \ldots , b_m)
$$
for any sequences
$(a_1 , \ldots , a_n)$, $(b_1 , \ldots , b_m)$ or real numbers.
In addition, for any function $g:\real^n \to \real$,
we let $g^*$ be the operator mapping
$C^{0, \infty}([0, T] \times \real^d)$
to
$C^{0, \infty}([0, T] \times \real^d)$
and defined by
$$
g^*(u)(t, x) := g\big(\partial_{\lambda^1}u(t,x),\ldots , \partial_{\lambda^n}u(t,x)\big), \qquad
 (t,x) \in [0,T]\times \real^d.
$$ 
In the sequel, we also let
$\partial_\lambda: = \partial_{z_1}^{\lambda_1} \cdots \partial_{z_n}^{\lambda_n}$,
and
$\partial_\mu : = \partial_{x_1}^{\mu_1} \cdots \partial_{x_d}^{\mu_d}$,
$\lambda = ( \lambda_1,\ldots , \lambda_n) \in \inte^n$,
$\mu = ( \mu_1,\ldots , \mu_d ) \in \inte^d$. 
\begin{definition}
 We let $\mathcal{C}$ denote the set of
 operators from ${\cal C}^{0,\infty} ([0,T]\times \real^d)$ to
 ${\cal C}^{0,\infty} ([0,T]\times \real^d)$,
 called \textit{codes}, and defined as
$$
 \mathcal{C} := \left\{
           {\rm Id}, \
           (a \partial_\lambda f)^*, \
           \partial_{\mu}, \ : \ 
           \lambda \in \mathbb{N}^n, \
           \mu \in \mathbb{N}^d, \
           a \in \real
                      \right\},
$$
where
${\rm Id}$
denotes the identity on ${\cal C}^{0,\infty} ([0,T]\times \real^d)$.
\end{definition}
\noindent
For example, for
$\nu \in \mathbb{N}^n$,
$\mu \in \mathbb{N}^d$,
$a \in \real$ and $k \in \mathbb{N}$
we have
\begin{equation*}
    c(u)(T, x) =
    \begin{cases}
      \displaystyle
      \phi(x),
                & \text{if } c = \rm Id,
      \medskip
      \\
      a \partial_\nu f\big(
      \partial_{\lambda^1}\phi(x)
, \ldots ,
      \partial_{\lambda^m}\phi(x)
      \big),
                & \text{if } c = (a \partial_\nu f)^*,
      \medskip
      \\
      \partial_\mu \phi(x),
                & \text{if } c = \partial_\mu.
    \end{cases}
\end{equation*}
 \noindent
 The mechanism $\mathcal{M}$ is then defined as a mapping on $\mathcal{C}$ by
$\mathcal{M} ( {\rm Id} ) := \{ f^* \}$, and
 \begin{align*}
 &
      {\cal M}(g^*)
 :=
       \bigcup_{1 \leq p \leq n \atop
       \lambda^p = 0}
       \big\{
       \big(
       f^*
       ,
       (\partial_{\bm{1}_p} g)^*
       \big)
       \big\}
       \\
    & 
 \bigcup_{
    \footnotesize \substack{1 \le p \le n, \
               1 \le s \le \abs{\lambda^p} \\
               1 \le \nu_1 + \cdots + \nu_n \le \abs{\lambda^p} \\
               1 \le \abs{k_1}, \dots, \abs{k_s}, \\
               0 \prec l^1 \prec \cdots \prec l^s \\
               k^i_1 + \cdots + k^i_s = \nu_i, \
               i = 1, \dots, n \\
               \abs{k_1}l_j^1 + \cdots + \abs{k_s}l_j^s = \lambda^p_j, \
               j = 1, \dots, d
                }}
\left\{
\left((\partial_{\bm{1}_p} g)^*,
\frac{\prod\limits_{i = 1}^d \lambda^p_i!  (\partial_{\nu} f)^*}
     {\prod\limits_{\footnotesize \substack{1 \le r \le s \\
                     1 \le q \le n}}
      k_r^q! \left(l_1^r! \cdots l_d^r!\right)^{k_r^q}
     }
\right)
\bigcup_{\footnotesize \substack{1 \le r \le s \\
                     1 \le q \le n \\
                    }}
\big(
\underbrace{
\partial_{l^r + \lambda^q} ,
\ldots
,
\partial_{l^r + \lambda^q} 
}
_{k_r^q~{\rm times}}
\big)
\right\}
\\
& \qquad \qquad \quad
\bigcup
\bigcup_{
     \footnotesize \substack{i, j = 1, \dots, n \\
               k = 1, \dots, d
                }}
\left\{
\left(-\frac{1}{2} (\partial_{\bm{1}_i + \bm{1}_j} g)^*,
    \partial_{\lambda^i + \bm{1}_k} ,
    \partial_{\lambda^j + \bm{1}_k} 
\right)\right\},
\qquad g^* \in \mathcal{C},
\end{align*}
 and
\begin{align*}
 &
  {\cal M}(\partial_{\mu})
\\
 &:=
 \bigcup_{
     \footnotesize \substack{1 \le s \le \abs{\mu}, \
               1 \le \nu_1 + \cdots + \nu_n \le \abs{\mu} \\
               1 \le \abs{k_1}, \dots, \abs{k_s}, \
               0 \prec l^1 \prec \cdots \prec l^s \\
               k^i_1 + \cdots + k^i_s = \nu_i, \
               i = 1, \dots, n \\
               \abs{k_1}l_j^1 + \cdots + \abs{k_s}l_j^s = \mu_j, \
               j = 1, \dots, d
                }}
\left\{
\left(
\frac{\prod\limits_{i = 1}^d \mu_i! }
     {\prod\limits_{\footnotesize \substack{1 \le r \le s \\
                     1 \le q \le n}}
      k_r^q! \left(l_1^r! \cdots l_d^r!\right)^{k_r^q}
     }
     (\partial_{\nu} f)^*
     \right)
     \bigcup_{\footnotesize \substack{1 \le r \le s \\
                     1 \le q \le n
                    }}
     \big(
     \underbrace{
       \partial_{l^r + \lambda^q} u
       ,
       \ldots
       ,
       \partial_{l^r + \lambda^q} u}_{k_r^q~{\rm times}}
     \big)
\right\},
\end{align*}
$\mu \in \mathbb{N}^d$.
 Given $\rho: \real_+ \to (0, \infty)$ 
 a probability density function (PDF) on $\real_+$
 with tail distribution function $\widebar{F}$
and $\mathcal{N}(0, \sigma^2)$ a $d$-dimensional
independent centered normal distribution with variance $\sigma^2$,
we consider the functional $\mathcal{H}({t, x, c})$
constructed in Algorithm~\ref{alg:coding tree}
 along a random coded tree started at
 $(t,x,c)\in [0,T]\times \real^d \times {\cal C}$,
 using independent random samples on a probability space $\Omega$.
\begin{algorithm}[H]
\caption{Coding tree algorithm TREE$(t, x, c)$}\label{alg:coding tree}
    \begin{algorithmic}
        \Require $t \in [0, T]$, $x \in \real^d$, $c \in \mathcal{C}$
        \Ensure $\mathcal{H}({t, x, c}) \in \real$
        \State $\mathcal{H}({t, x, c}) \gets 1$
        \State $\tau \gets$ a random variable drawn from the distribution of $\rho$
        \If{$t + \tau > T$}
            \State $W \gets$ a random vector drawn from $\mathcal{N}(0, T - t)$
            \State $\mathcal{H}({t, x, c}) \gets \mathcal{H}({t, x, c}) \times c(u)(T, x + W)
            /  \widebar{F}(T-t)$
        \Else
            \State $q \gets$ the size of the mechanism set $\mathcal{M}(c)$
            \State $I \gets$ a random element drawn uniformly from $\mathcal{M}(c)$
            \State $\mathcal{H}({t, x, c}) \gets \mathcal{H}({t, x, c}) \times q / \rho(\tau)$
            \ForAll{$cc \in I$}
                \State $W \gets$ a random vector drawn from $\mathcal{N}(0, \tau)$
                \State $\mathcal{H}({t, x, c}) \gets \mathcal{H}({t, x, c}) \times \text{TREE}(t + \tau, x + W, cc)$
            \EndFor
        \EndIf
    \end{algorithmic}
\end{algorithm}
\noindent
 As in Theorem~1 in \cite{penent2022fully},
 the following Feynman-Kac type identity
\begin{equation}
    \label{eq:feynman kac}
    u(t, x) = \E [ {\cal H} ({t, x, \rm Id} ) ] 
\end{equation}
for the solution of \eqref{eq:main pde} holds under suitable integrability conditions
on ${\cal H} ({t, x, \rm Id} )$ and smoothness assumptions on the coefficients of
 \eqref{eq:main pde}, see the appendix for calculation details.
\section{Deep branching solver}
\label{sec:main ideas}
Instead of evaluating \eqref{eq:feynman kac} at a single point
$(t,x) \in [0,T]\times \real^d$, we use the $L^2$-minimality property of expectation to perform a functional estimation of $u(\cdot , \cdot )$
as $u(\cdot , \cdot ) = v^*(\cdot , \cdot )$
 on the support of a random vector $(\tau , X)$ on $[0,T] \times \real^d$
 such that $\mathcal{H} ({\tau, X, \rm Id}) \in L^2 (\Omega )$, where
\begin{equation}
\label{eq:L2 minimality}
    v^* = \argmin\limits_{\{ v \ \! : \ \! v(\tau, X) \in L^2\}} \E \left[
            \left(\mathcal{H}({\tau, X, \rm Id}) - v(\tau, X)\right)^2 \right].
\end{equation}
To evaluate \eqref{eq:feynman kac} on $[0, T] \times \domain$,
where $\domain$ is a bounded domain of $\real^d$,
we can choose $(\tau, X)$ to be a uniform random vector on $[0, T] \times \domain$.
Similarly, to evaluate \eqref{eq:feynman kac} on $\{0\} \times \domain$,
we may let $\tau \equiv 0$ and let $X$ be a uniform random vector on $\domain$.

\medskip

In order to implement the deep learning approximation,
we parametrize $v(\cdot , \cdot )$
 in the functional space described below.
Given $\sigma:\real \to \real$ an activation function
such as $\sigma_{\rm ReLU}(x) := \max (0,x)$,
$\sigma_{\tanh}(x) := \tanh(x)$ or $\sigma_{\rm Id}(x) := x$,
 we define the set of layer functions $\mathbb{L}^\sigma_{d_1,d_2}$ by
\begin{equation}
    \label{eq:one layer}
    \mathbb{L}^\sigma_{d_1,d_2} :=
    \bigl\{
        L:\real^{d_1} \to \real^{d_2} \ : \ L(x) = \sigma( Ax + b),
  \ x \in \real^{d_1}, \ A \in \real^{d_2 \times d_1}, \ b \in \real^{d_2}
    \bigr\},
\end{equation}
where $d_1 \geq 1$ is the input dimension,
$d_2 \geq 1$ is the output dimension,
and the activation function
$\sigma$ is applied component-wise to $Ax + b$.
Similarly, when the input and output dimensions are the same,
we define the set of residual layer functions
$\mathbb{L}^{\rho, \rm res}_d$ by
\begin{equation}
    \label{eq:one layer resnet}
    \mathbb{L}^{\sigma, \rm res}_d :=
    \bigl\{
        L:\real^d \to \real^d \ : \ L(x) = x + \sigma(Ax + b),
  \ x \in \real^d, \ A \in \real^{d \times d}, \ b \in \real^d
    \bigr\},
\end{equation}
see \cite{He16}.
Then, we denote by
\begin{equation*}
    \mathbb{NN}^{\sigma,l,m}_d :=
    \bigl\{
    L_l \circ \dots \circ L_0 : \real^d \to \real
    \ : \
        L_0 \in \mathbb{L}^{\sigma}_{d, m}, \ 
        L_l \in \mathbb{L}^{\sigma_{\rm Id}}_{m, 1}, \ 
        L_i \in \mathbb{L}^{\sigma, \rm res}_m, \ 
        1 \leq i < l
    \bigr\}
\end{equation*}
the set of feed-forward neural networks
with one output layer,
$l \geq 1$ hidden residual layers
each containing $m \geq 1$ neurons,
 where the activation functions of
the output and hidden layers
 are respectively
the identity function $\sigma_{\rm Id}$
and $\sigma$.
Any $v(\cdot; \theta) \in \mathbb{NN}^{\sigma,l,m}_d$
is fully determined by the sequence
\begin{equation*}
    \theta := \bigl( A_0, b_0, A_1, b_1, \dots, A_{l-1}, b_{l-1}, A_l, b_l \bigr)
\end{equation*}
of $\left( (d+1) m + (l - 1) (m+1) m + (m+1) \right)$ parameters.

\medskip

Since by the universal approximation theorem,
see e.g.\ Theorem~1 of \cite{hornik1991approximation},
$\bigcup\limits_{m = 1}^{\infty} \mathbb{NN}^{\sigma,l ,m}_d$
is dense in the $L^2$ functional space,
the optimization problem \eqref{eq:L2 minimality} can be approximated by
\begin{equation}
\label{eq:NN approximation}
    v^* \approx \argmin\limits_{v \in \mathbb{NN}^{\sigma,l ,m}_{d + 1}}
            \E \left[
                \left(\mathcal{H}({\tau, X, \rm Id}) - v(\tau, X)\right)^2 \right].
\end{equation}
By the law of large numbers,
\eqref{eq:NN approximation} can be further approximated by
\begin{equation}
\label{eq:monte carlo}
    v^* \approx \argmin\limits_{v \in \mathbb{NN}^{\sigma,l ,m}_{d + 1}}
        N^{-1} \sum\limits_{i = 1}^N
            \left(\mathcal{H}_i - v(\tau_i, X_i)\right)^2,
\end{equation}
where for all $i = 1, \dots, N$,
$(\tau_i, X_i)$ is drawn independently from the distribution of $(\tau, X)$
and $\mathcal{H}_i$ is drawn from $\mathcal{H}_{\tau_i, X_i, \rm Id}$
using Algorithm~\ref{alg:coding tree}.
However,
the approximation \eqref{eq:monte carlo} may perform poorly
when the variance of $\mathcal{H}_i$ is too high.
To solve this issue, we use the expression 
\begin{equation}
\label{eq:monte carlo more samples}
    v^* \approx \argmin\limits_{v \in \mathbb{NN}^{\sigma,l ,m}_{d + 1}}
        N^{-1} \sum\limits_{i = 1}^N
            \left(M^{-1} \sum\limits_{j = 1}^M \mathcal{H}_{i,j} - v(\tau_i, X_i)\right)^2,
\end{equation}
 where for $j = 1, \dots, M$, $\mathcal{H}_{i, j}$ is drawn independently
 from $\mathcal{H}_{\tau_i, X_i, \rm Id}$ using Algorithm~\ref{alg:coding tree}.

\medskip

Finally, the deep branching method
using the gradient descent method
to solve the optimization in \eqref{eq:monte carlo more samples}
is summarized in Algorithm~\ref{alg:deep}.
\begin{algorithm}[H]
\caption{Deep branching method}\label{alg:deep}
    \begin{algorithmic}
        \Require The learning rate $\eta$ and the number of epochs $P$
        \Ensure $v(\cdot, \cdot; \theta) \in \mathbb{NN}^{\sigma,l,m}_{d+1}$
        \State $(\tau_i, X_i)_{1 \le i \le N} \gets$ random vectors
                    drawn from the distribution of $(\tau, X)$
        \State $(\mathcal{H}_{i,j})_{\footnotesize \substack{1 \le i \le N \\ 1 \le j \le M}} \gets$
                    random variables
                    generated by TREE$(\tau_i, X_i, \rm Id)$ in
                    Algorithm~\ref{alg:coding tree}
        \smallskip
        \State Initialize $\theta$
        \For{$i \gets 1, \dots, P$}
            \State $L \gets N^{-1} \sum\limits_{i = 1}^N
                \bigg( M^{-1} \sum\limits_{j = 1}^M \mathcal{H}_{i,j}
                    - v(\tau_i, X_i; \theta)\bigg)^2$
            \State $\theta \gets \theta - \eta \nabla_\theta L$
        \EndFor
    \end{algorithmic}
\end{algorithm}
\begin{remark}
    \label{remark:tricks}
 In the implementation of Algorithm~\ref{alg:deep},
 we perform the following additional steps:
    \begin{enumerate}[i)]
        \item $\eta \gets \eta / 10$ after every $\floor{P/3}$ steps.
        \item Instead of using $\eta$ to update $\theta$ directly,
                Adam algorithm is used to update $\theta$,
                see \cite{kingma2014adam}.
        \item $\sigma_{\tanh}$ is used because
                the target PDE solution \eqref{eq:main pde} is smooth.
        \item A batch normalization layer
                is added after the activation function in
                \eqref{eq:one layer}-\eqref{eq:one layer resnet}
                when $\sigma \neq \sigma_{\rm Id}$,
                see \cite{ioffe2015batch}.
        \item $\rho$ is chosen to be the PDF of exponential distribution
                with rate $- ( \log 0.95) / T$.
        \item Given $x_{\rm min} < x_{\rm max}$ and
                $x_{\rm mid} = (x_{\rm min} + x_{\rm max})/2$,
          we take
          $$
          \domain := [x_{\rm min}, x_{\rm max}]
          \times \{ x_{\rm mid} \}\times \dots \times \{ x_{\rm mid} \},
          $$
                and we let $(\tau, X)$ be the uniform random vector on $\{ 0 \} \times \domain$.
    \end{enumerate}
\end{remark}

\section{Numerical examples}
\label{sec:numerical examples}
The numerical examples below are run in {\sc Python} using {\sc PyTorch}
with the default initialization scheme for $\theta$, and the default values
$N=1000$, $P=3000$, $\eta=0.01$, $l=6$, $m=20$.
Except if otherwise stated, runtimes are expressed in minutes and the examples have been run on Google Colab with a Tesla P100 GPU. 

\medskip

For comparisons with the deep BSDE and deep Galerkin methods,
we select the configurations such that all methods have comparable or similar runtimes.
For the deep BSDE method of \cite{han2018solving,beck},
the time discretization of
$(0, {T} / {5}, {2T} / {5}, {3T} / {5}, {4T} / {5}, T )$
and $1000$ (resp. $100,000$) number of samples
are used in the case of $d=1$ (resp. $d>1$).

For the deep Galerkin method of \cite{sirignano2018dgm},
$10,000$ samples are respectively generated on
$\{ 0 \} \times [x_{\rm min}, x_{\rm max}]^d$,
$(0, T) \times [x_{\rm min}, x_{\rm max}]^d$,
and $\{ T \} \times [x_{\rm min}, x_{\rm max}]^d$.
In our experiment, such generation works better
than generating $10,000$ samples respectively on
$\{ 0 \} \times \domain$,
$(0, T) \times \domain$,
and $\{ T \} \times \domain$.
In addition, we found that batch normalization and learning rate decay
in Remark~\ref{remark:tricks}
do not work well with deep Galerkin method,
hence they are not used in the simulation below for the deep Galerkin method.
The learning rate for the deep Galerkin method is fixed to be $\eta = 0.001$
throughout the training.

The analysis of error is performed on the grid of
$\widetilde{\domain} = (0, x_{\rm min} + i \Delta_x, x_{\rm mid}, \dots, x_{\rm mid})_{0 \le i \le 100}$,
where $\Delta_x = (x_{\rm max} - x_{\rm min}) / 100$.
In each of the $10$ independent runs,
the statistics of the runtime (in seconds) and
the $L^p$ error
$100^{-1} \sum\limits_{x \in \widetilde{\domain}}
\abs{{\rm true}(x) - {\rm predicted}(x)}^p$
are recorded.
In multidimensional examples with $d\geq 2$,
every figure is plotted as a function of $x_1$
on the horizontal axis, after setting $(x_2,\ldots , x_d)=(0,\ldots ,0)$. 

\subsubsection*{a) Allen-Cahn equation}
\label{subsec:example 1}
\begin{figure}[H]
\centering
\includegraphics[width=0.55\textwidth]{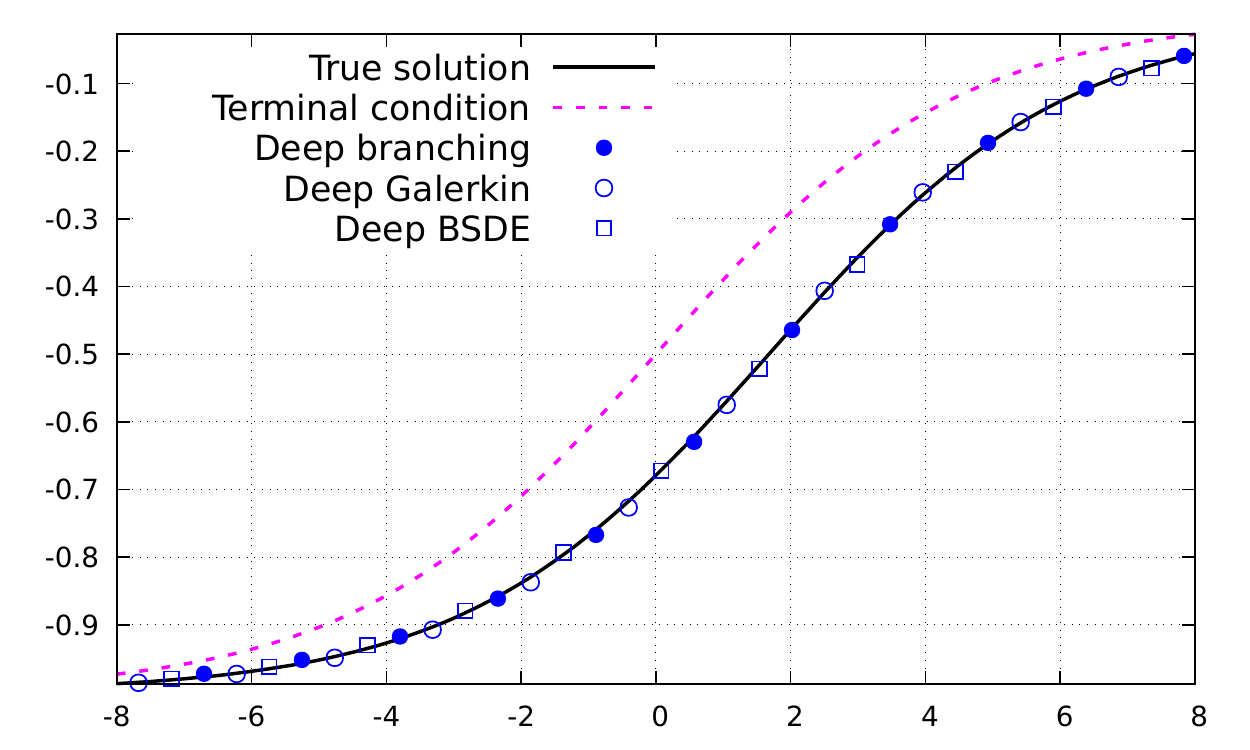}
\vskip-0.1cm
\caption{Comparison of deep learning methods for \eqref{eq:example 1} with $d = 5$ and $T=0.5$.}
\label{fig:example 1}
\end{figure}

Consider the equation
\begin{equation}
\label{eq:example 1}
\partial_t u(t,x) + \frac{1}{2} \Delta u(t,x) + u(t,x) - u^3(t,x) = 0,
\end{equation}
which admits the traveling wave solution
$$u(t,x) = -\frac{1}{2} - \frac{1}{2}
\tanh \left( \frac{3}{4} (T-t) - \sum\limits_{i=1}^d \frac{x_i}{2\sqrt{d}} \right),
\qquad
(t,x) \in [0,T]\times \real^d.
$$
Table~\ref{table:example 1} summarizes the results of $10$ independent runs,
with $M = 100,000$, $T = 0.5$, $x_{\rm min} = -8$, and $x_{\rm max} = 8$.
\begin{table}[H]
    \centering
    \resizebox{\textwidth}{!}{\begin{tabular}{|l|c||c|c||c|c||c|}
        \hline
\multicolumn{1}{|c|}{Method}
         & $d$ & Mean $L^1$-error & Stdev $L^1$-error &
        Mean $L^2$-error & Stdev $L^2$-error & Mean Runtime \\
        \hline
        Deep branching                        & 1 & 1.32E-03 & 1.05E-04 & 4.04E-06 & 7.32E-07 & 28m \\
         \cline{1-1}\cline{3-7}
        Deep BSDE \cite{han2018solving}       & 1 & 4.60E-03 & 9.82E-04 & 4.08E-05 & 2.18E-05 & 101m \\
         \cline{1-1}\cline{3-7}
        Deep Galerkin \cite{sirignano2018dgm} & 1 & 1.40E-03 & 1.83E-03 & 6.39E-06 & 1.57E-05 & 53m \\
        \hline
        Deep branching                        & 5 & 3.63E-03 & 1.57E-04 & 2.09E-05 & 1.19E-06 & 110m \\
         \cline{1-1}\cline{3-7}
        Deep BSDE \cite{han2018solving}       & 5 & 4.71E-03 & 4.23E-04 & 3.51E-05 & 8.19E-06 & 170m \\
         \cline{1-1}\cline{3-7}
        Deep Galerkin \cite{sirignano2018dgm} & 5 & 6.83E-03 & 6.17E-03 & 1.36E-04 & 2.77E-04 & 134m \\
        \hline
\end{tabular}}
		\caption{Summary of numerical results for \eqref{eq:example 1}.}
        \label{table:example 1}
\end{table}
\vskip-0.3cm
\noindent
 We check in Table~\ref{table:example 1} and Figure~\ref{fig:example 1} 
 that all three algorithms show a similar accuracy for the numerical solution
 of the Allen-Cahn equation, while the deep branching method appears more
 stable.

\vskip-0.1cm

\vskip-0.2cm

\begin{figure}[H]
  \centering
 \begin{subfigure}[b]{0.49\textwidth}
\includegraphics[width=\textwidth]{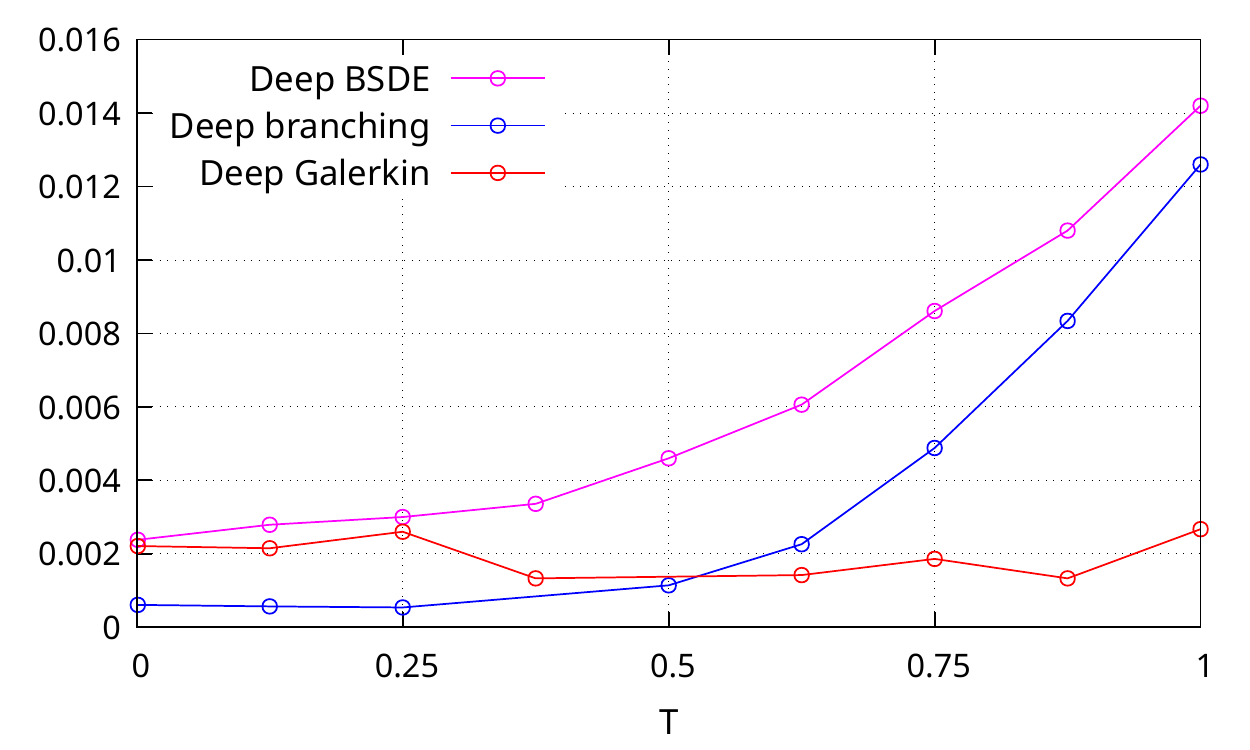} 
\caption{Dimension $d=1$.}
 \end{subfigure}
  \begin{subfigure}[b]{0.49\textwidth}
\includegraphics[width=\textwidth]{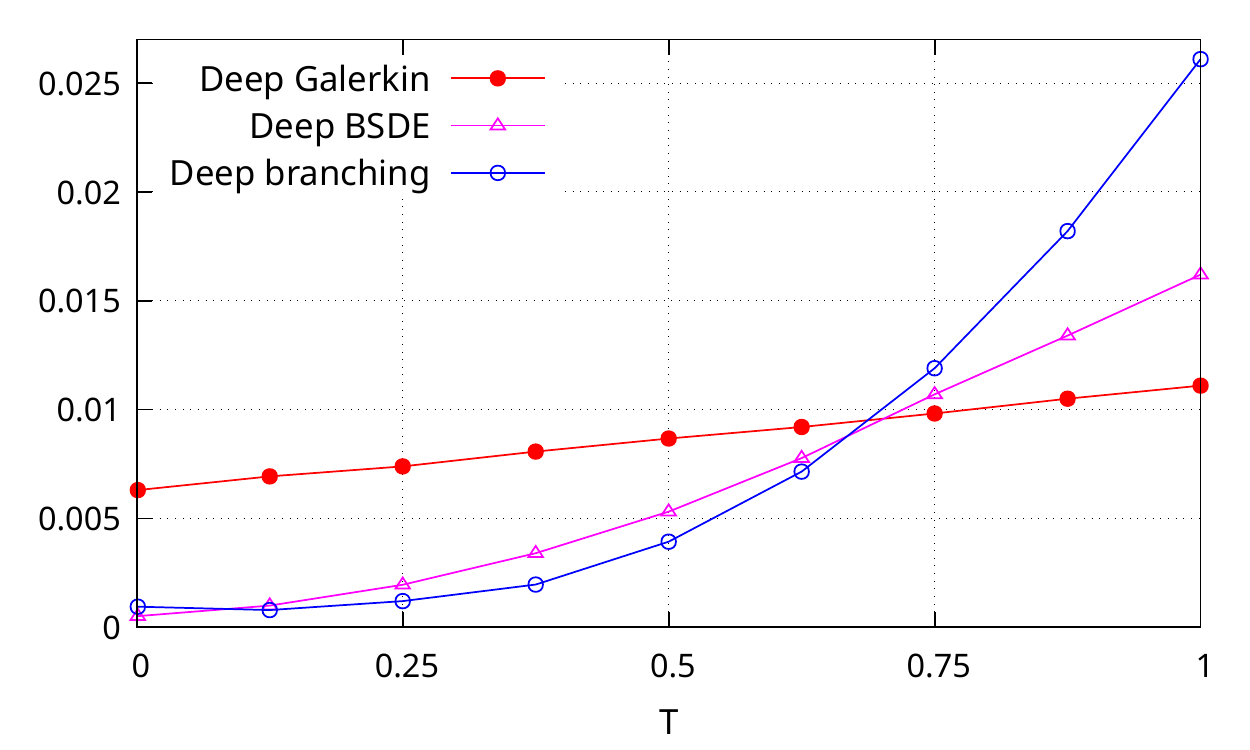} 
\caption{Dimension $d=10$.}
 \end{subfigure}
  \caption{$L^1$ error graphs for \eqref{eq:example 1} as functions of time $T$.} 
\label{figexpl}
\end{figure}

\begin{figure}[H]
\centering
\includegraphics[width=0.55\textwidth]{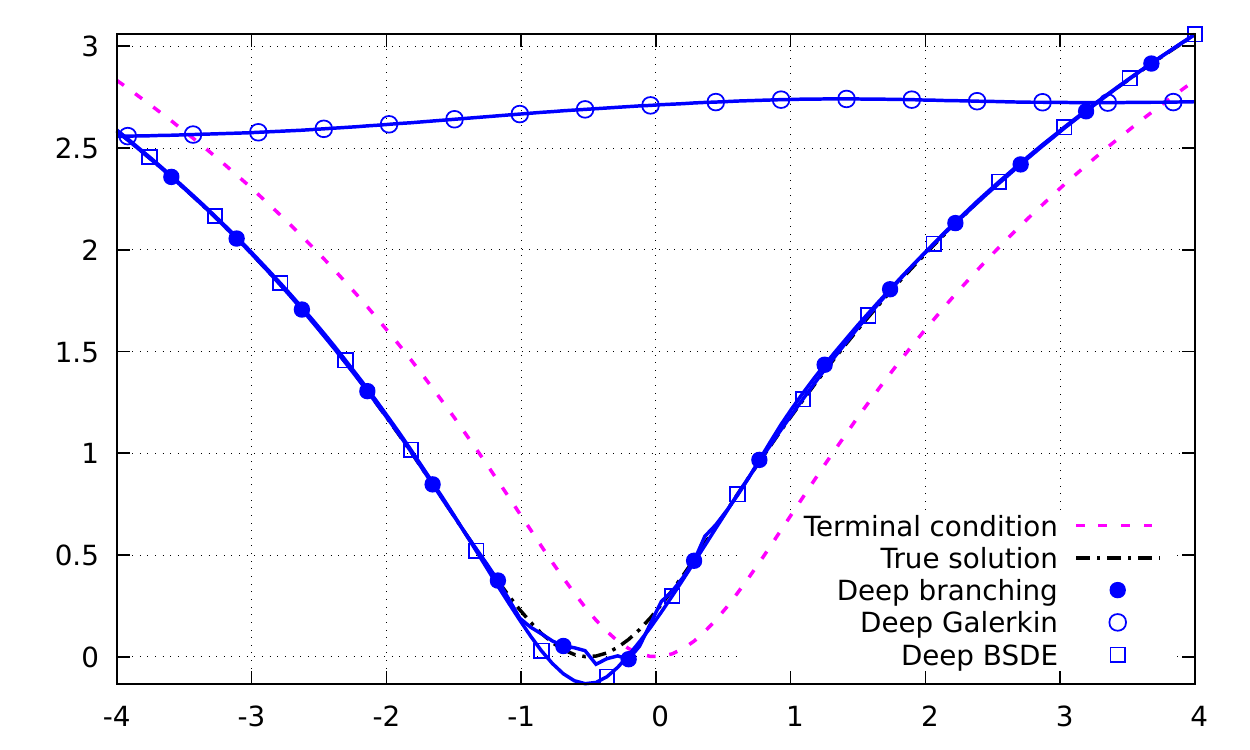}
\vskip-0.1cm
\caption{Comparison of deep learning methods
  for \eqref{eq:example 2} with $d = 5$ and $T=0.05$.}
\label{fig:example 2}
\end{figure}

\noindent
{Figure~\ref{figexpl} compares the $L^1$ errors 
  of deep learning methods, showing that although the 
  deep branching method has an explosive behavior,
  under comparable runtimes
  it can perform better than the deep Galerkin and
  deep BSDE methods in small time, in both dimensions $d=1$
  and $10$.
  Figure~\ref{figexpl} and Table~\ref{table2}
  have been run on a RTX A4000 GPU.  
  }

\vskip-0.4cm

\noindent
Table~\ref{table2} ensures that the experiments of
Figure~\ref{figexpl} are performed within comparable
runtimes. 
    
\begin{table}[H]
    \centering
            {\begin{tabular}{|l||c|c|}
        \hline
        \multicolumn{1}{|c||}{Method} & $d=1$ & $d=10$ \\
        \hline
        Deep branching                         & 38m & 279m \\
        \hline
        Deep Galerkin \cite{han2018solving}        & 69m & 254m \\
        \hline
        Deep BSDE \cite{sirignano2018dgm}  & 87m & 237m \\
        \hline
\end{tabular}}
	\caption{Average runtimes in minutes for Figure~\ref{figexpl}.
        }
\label{table2} 
\end{table}

\vskip-0.3cm

\noindent

\subsubsection*{b) Exponential nonlinearity}
\label{subsec:example 2}
Consider the equation
\begin{equation}
\label{eq:example 2}
\partial_t u(t,x) + \frac{\alpha}{d} \sum\limits_{i=1}^d \partial_{x_i} u(t,x)
+ \frac{1}{2} \Delta u(t,x)
+ \re^{-u(t,x)} ( 1 - 2 \re^{-u(t,x)} ) d = 0,
\end{equation}
which admits the traveling wave solution
$$
u(t,x) = \log \left( 1 + \left( \sum\limits_{i=1}^d x_i + \alpha (T-t)\right)^2 \right),
\qquad (t,x) \in [0,T]\times \real^d.
$$
Table~\ref{table:example 2} summarizes the results of $10$ independent runs,
with $M = 30,000$ (resp. $M = 3,000$) in dimension $d = 1$ (resp. $d = 5$),
$\alpha = 10$, $T = 0.05$, $x_{\rm min} = -4$, and $x_{\rm max} = 4$.

\begin{table}[H]
    \centering
    \resizebox{\textwidth}{!}{\begin{tabular}{|l|c||c|c||c|c||c|}
        \hline
        \multicolumn{1}{|c|}{Method} & $d$ & Mean $L^1$-error & Stdev $L^1$-error &
        Mean $L^2$-error & Stdev $L^2$-error & Mean Runtime \\
        \hline
        Deep branching                        & 1 & 1.17E-02 & 1.36E-03 & 4.57E-04 & 1.32E-04 & 42m \\
        \cline{1-1}\cline{3-7}
        Deep BSDE \cite{han2018solving}       & 1 & 1.39E-02 & 2.26E-03 & 3.56E-04 & 1.03E-04 & 101m \\
        \cline{1-1}\cline{3-7}
        Deep Galerkin \cite{sirignano2018dgm} & 1 & 2.53E-02 & 2.12E-02 & 1.72E-03 & 3.02E-03 & 61m \\
        \hline
        Deep branching                        & 5 & 2.63E-02 & 4.53E-03 & 2.69E-03 & 1.08E-03 & 146m \\
        \cline{1-1}\cline{3-7}
Deep BSDE \cite{han2018solving}       & 5 & 1.88E-02 & 4.57E-04 & 1.36E-03 & 9.86E-05 & 119m \\
        \cline{1-1}\cline{3-7}
        Deep Galerkin \cite{sirignano2018dgm} & 5 & 1.32E+00 & 7.78E-01 & 3.26E+00 & 2.54E+00 & 154m \\
        \hline
\end{tabular}}
		\caption{Summary of numerical results for \eqref{eq:example 2}.}
        \label{table:example 2}
\end{table}
\vskip-0.3cm
\noindent 
In the case of exponential nonlinearity, our method appears
significantly more accurate than the deep Galerkin
method, and performs comparably to the deep BSDE method
in dimension $d=5$. 

\vskip-0.1cm

\subsubsection*{c) Burgers equation} 
\begin{figure}[H]
\centering
\hskip-0.35cm
\begin{subfigure}[b]{0.52\textwidth}
\includegraphics[width=\textwidth]{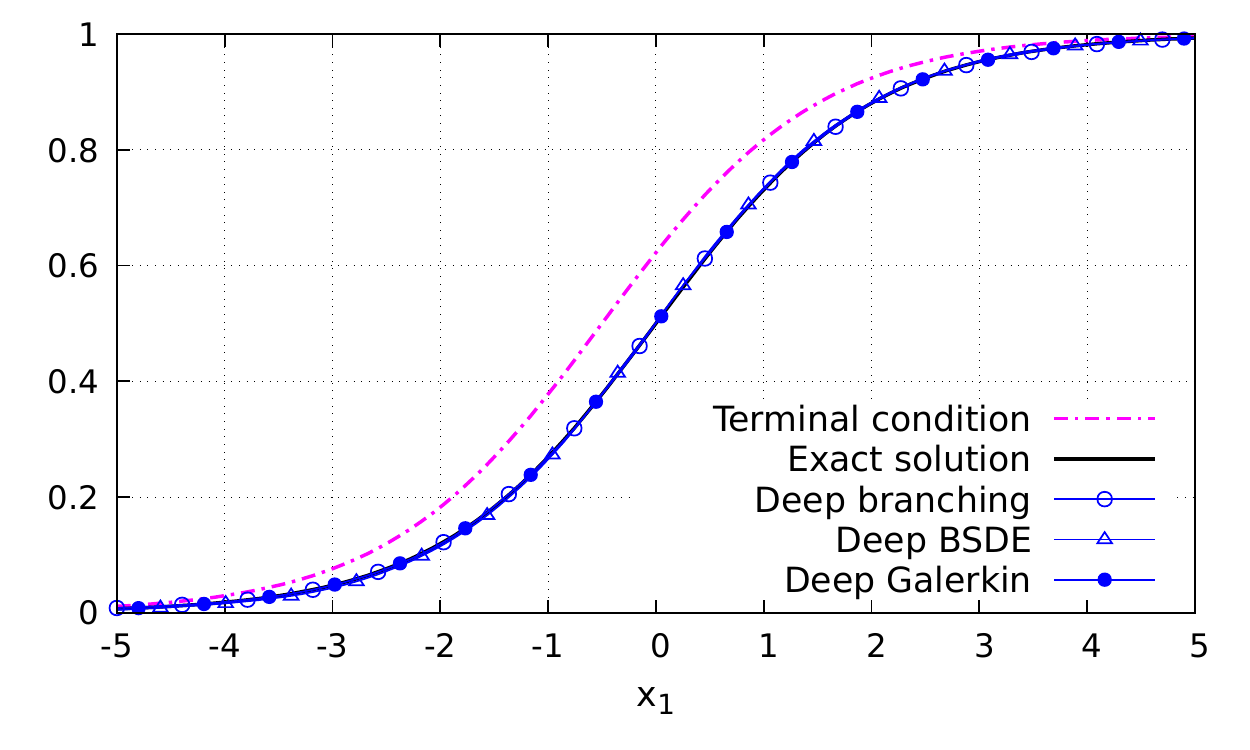} 
\caption{Dimension $d=1$ with $T=0.5$.}
 \end{subfigure}
\hskip-0.5cm
 \begin{subfigure}[b]{0.52\textwidth}
\includegraphics[width=\textwidth]{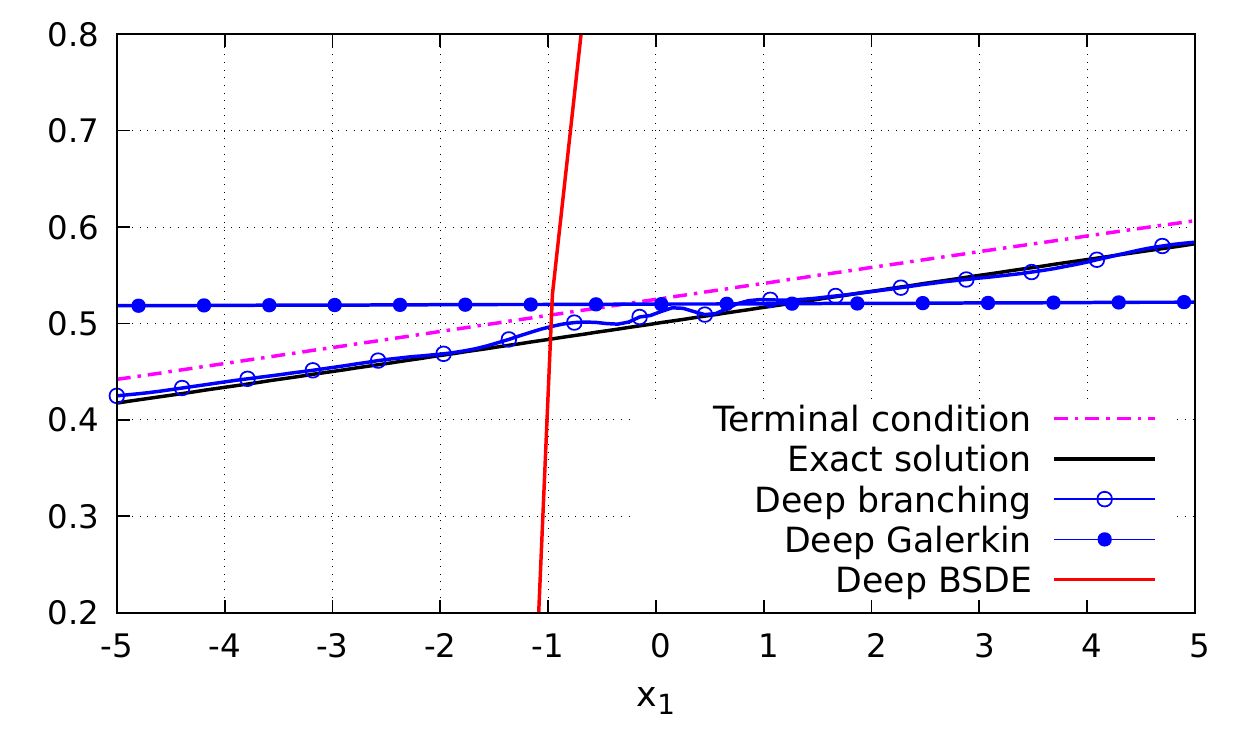} 
\caption{Dimension $d=15$ with $T=0.1$.}
 \end{subfigure}
   \caption{Numerical solution of \eqref{burgers0-00} and comparison to
   \eqref{tr2-00} with $\nu=d^2$.}
\label{fig1-00}
\end{figure}

{
\noindent
 Next, we consider the multidimensional Burgers equation 
 \label{p-13}
 \begin{equation}
    \label{burgers0-00} 
    \partial_t u(t,x) + \frac{d^2}{2} \Delta u(t,x)
    +
    \left( u (t,x) - \frac{2+d}{2d} \right) 
    \left(
    d \sum_{k=1}^d \partial_{x_k} u(t,x)
    \right)
    = 0,
\end{equation} 
 with traveling wave solution 
  \begin{equation}
   \label{tr2-00} 
   u(t, x) = \frac{\exp \big( t + d^{-1} \sum_{i=1}^d x_i
   \big)}{1 + \exp \big( t + d^{-1} \sum_{i=1}^d x_i
   \big)
   },
   \qquad
   x=(x_1,\ldots , x_d) \in \real^d, \ t\in [0,T],  
\end{equation} 
 see \S~4.5 of \cite{han2018solvingarxiv},
 and \S~4.2 of \cite{chassagneux}. 
  Figure~\ref{fig1-00} presents
  estimates of the solution of the Burgers equation \eqref{burgers0-00}
  with solution \eqref{tr2-00} in dimensions $d=5$ and $d=20$,
  with comparisons to the outputs 
  of the deep Galerkin method \cite{sirignano2018dgm}
  and of the deep BSDE method \cite{han2018solving}.}

\vskip-0.3cm

\begin{figure}[H]
  \centering
 \begin{subfigure}[b]{0.49\textwidth}
\includegraphics[width=\textwidth]{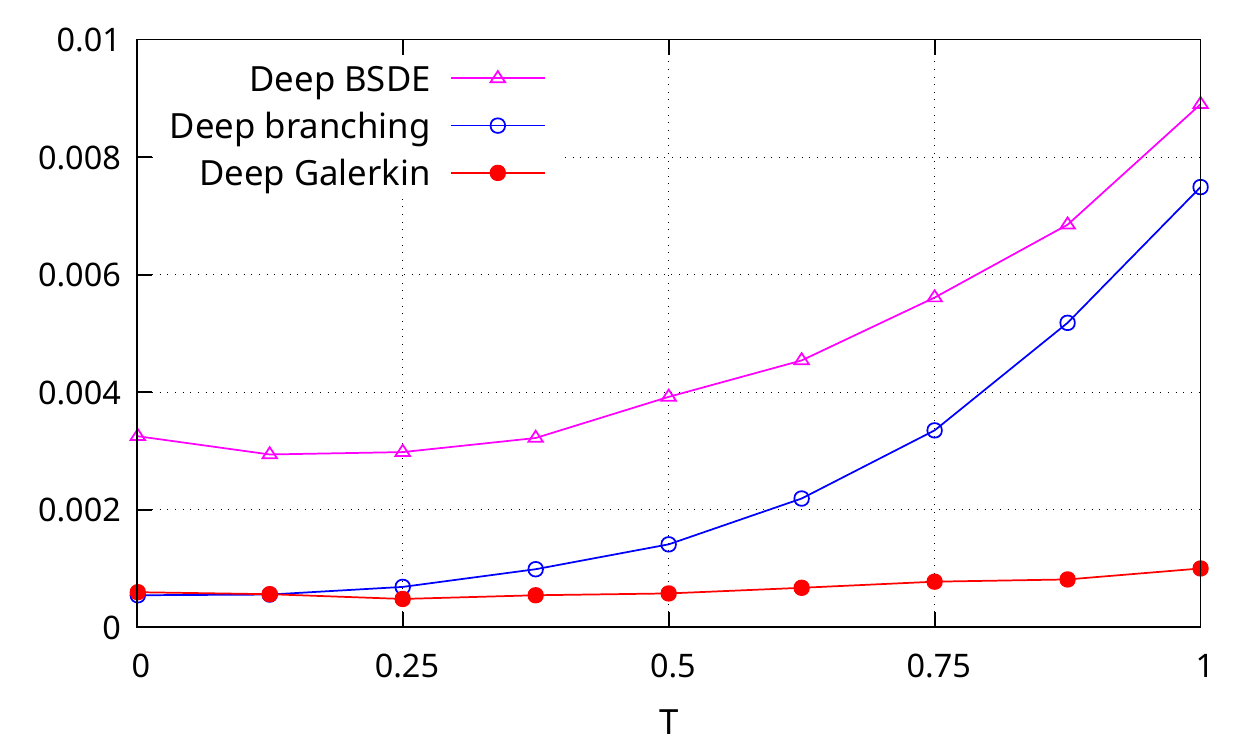} 
\caption{Dimension $d=1$.}
 \end{subfigure}
  \begin{subfigure}[b]{0.49\textwidth}
\includegraphics[width=\textwidth]{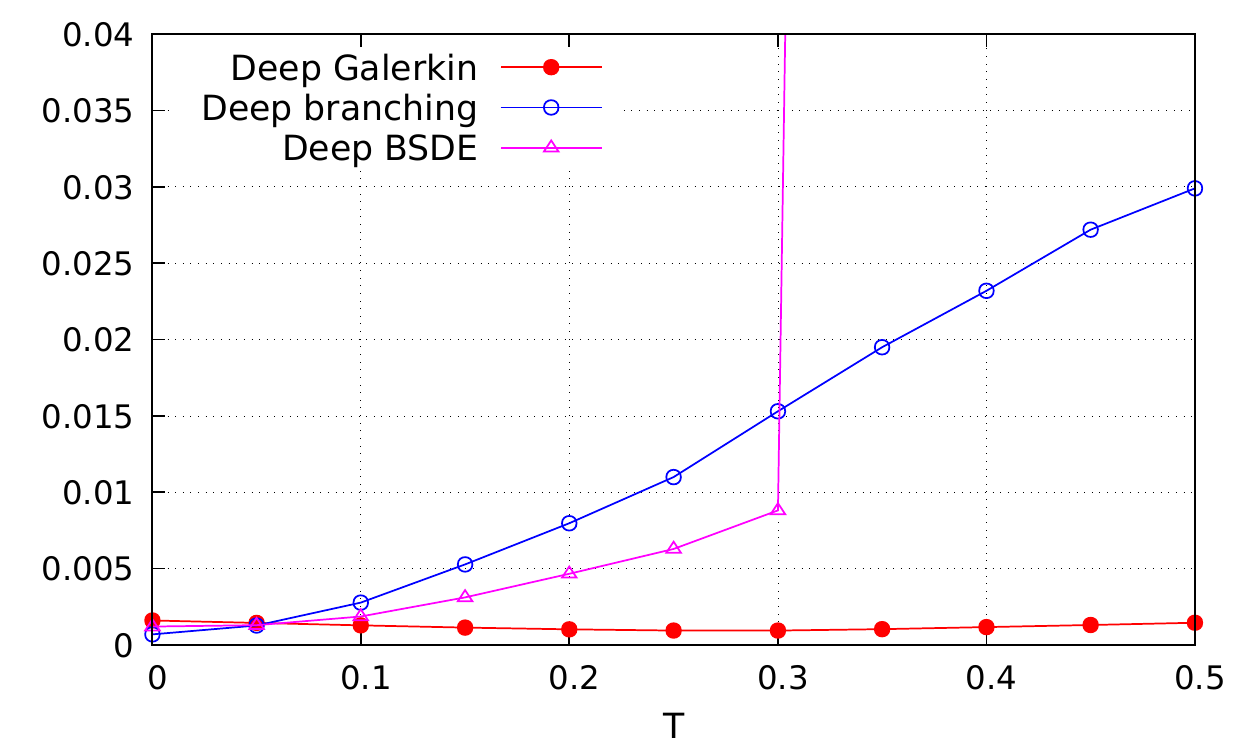} 
\caption{Dimension $d=5$.}
 \end{subfigure}
 \begin{subfigure}[b]{0.49\textwidth}
\includegraphics[width=\textwidth]{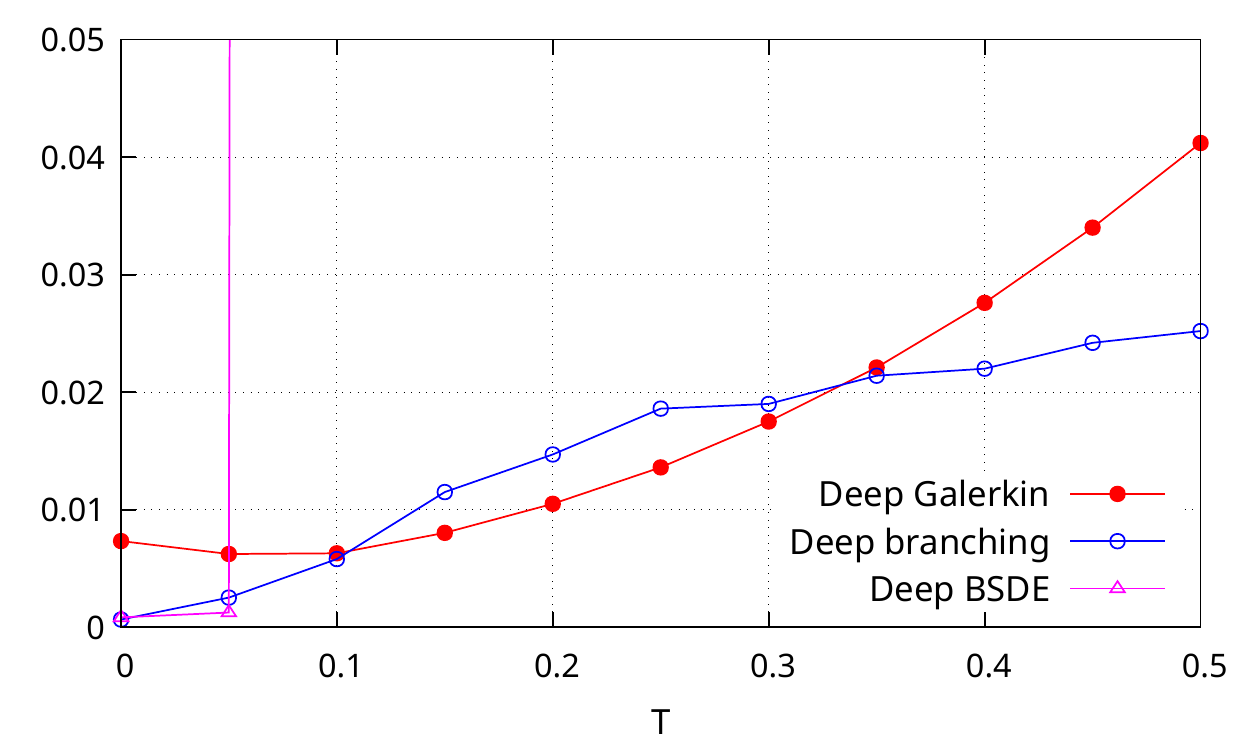} 
\caption{Dimension $d=10$.}
 \end{subfigure}
  \begin{subfigure}[b]{0.49\textwidth}
\includegraphics[width=\textwidth]{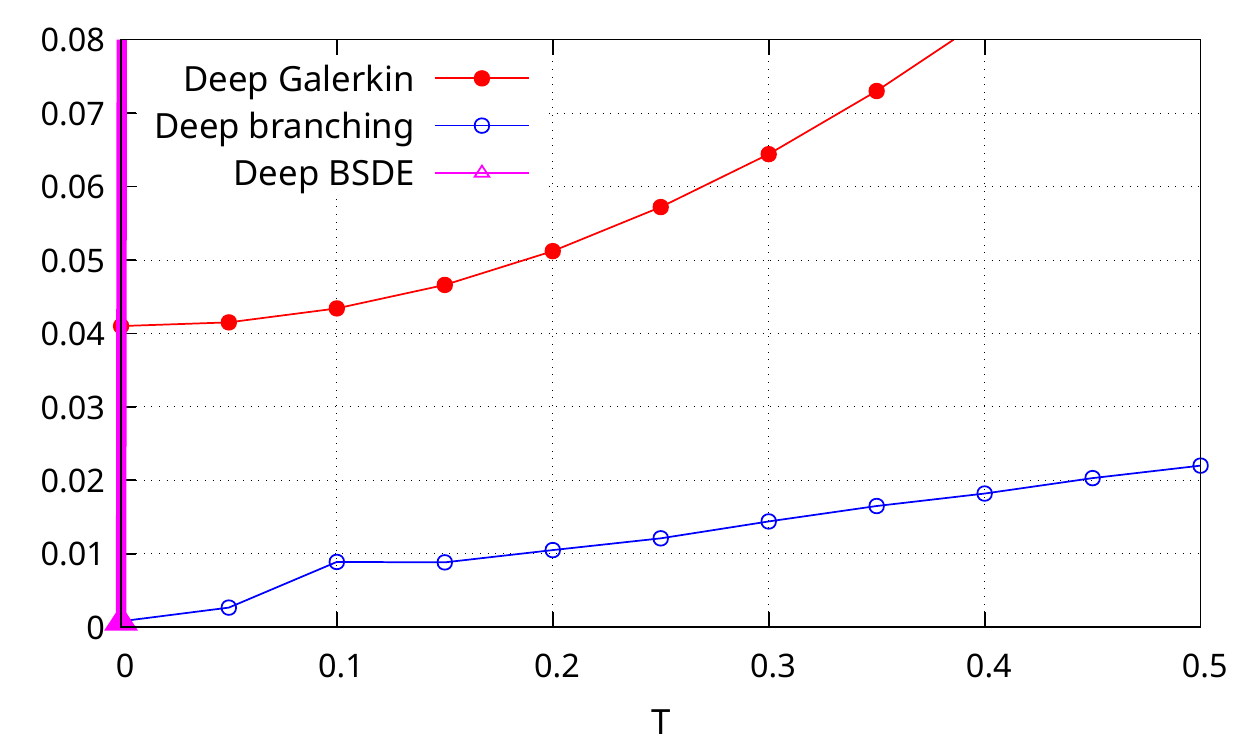} 
\caption{Dimension $d=15$.}
 \end{subfigure}
  \caption{$L^1$ error graphs for the solution of \eqref{burgers0-00} as functions of time $T$.} 
\label{figexpl-2}
\end{figure}

\begin{figure}[H]
\centering
\includegraphics[width=0.55\textwidth]{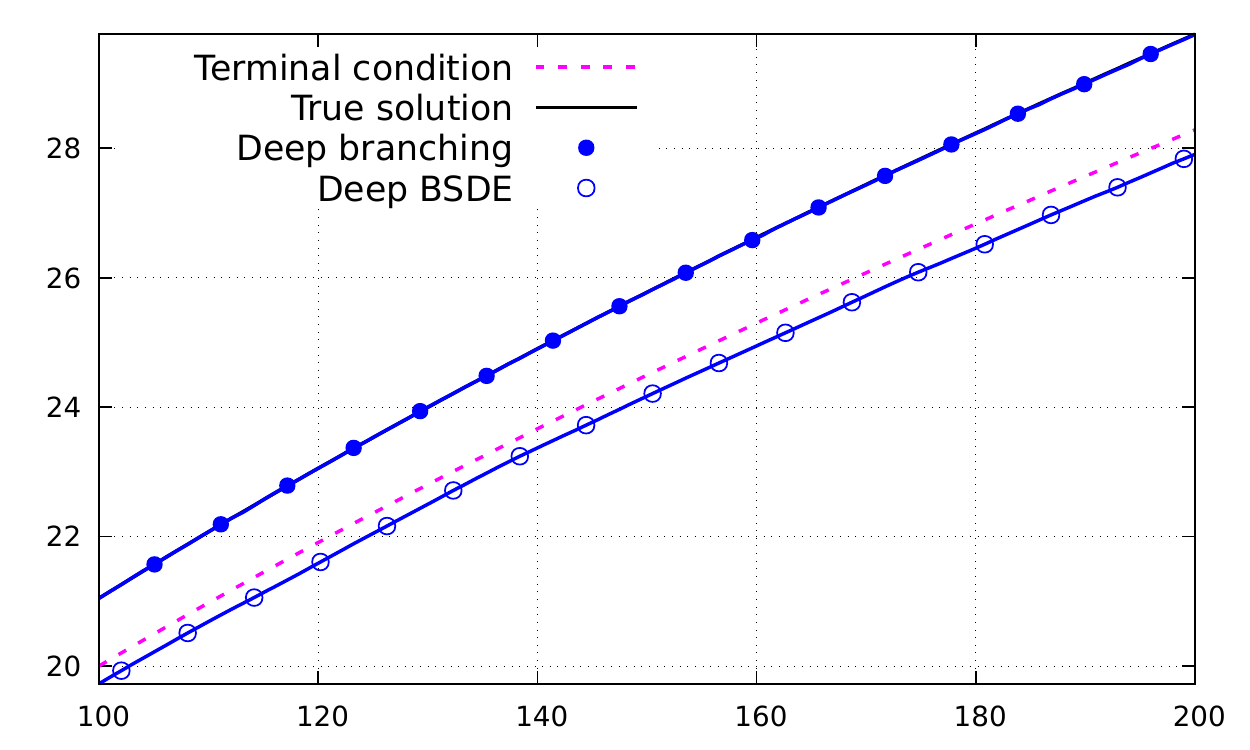}
\vskip-0.1cm
\caption{Deep branching {\em vs} deep BSDE method for
  \eqref{eq:example 7} with $d = 1$ and $T=0.1$.}
\label{fig:example 7}
\end{figure}

\noindent
{
 We note in Figure~\ref{fig1-00}$-b)$ that the deep branching method
 is more stable than the deep Galerkin and deep BSDE methods
 in dimension $d=15$. In particular, the deep BSDE estimate explodes
 under comparable runtimes, as shown in Figure~\ref{figexpl-2}.
 Figure~\ref{figexpl-2} and Table~\ref{table4} 
 have been run on a RTX A4000 GPU.
}

\vskip-0.4cm

\noindent
Table~\ref{table4} ensures that the experiments of
Figure~\ref{figexpl-2} are performed within comparable
runtimes. 
    
\begin{table}[H]
    \centering
            {\begin{tabular}{|l||c|c|c|c|}
        \hline
        \multicolumn{1}{|c||}{Method} & $d=1$ & $d=5$ & $d=10$ & $d=15$ \\
        \hline
        Deep branching                         & 61m & 72m & 117m & 169m \\
        \hline
        Deep Galerkin \cite{han2018solving}        & 77m & 184m & 331m & 480m \\
        \hline
        Deep BSDE \cite{sirignano2018dgm}  & 144m & 146m & 145m & 124m \\
        \hline
\end{tabular}}
	\caption{Average runtimes in minutes for Figure~\ref{figexpl-2}.
        }
\label{table4} 
\end{table}

\vskip-0.3cm

\noindent

\subsubsection*{d) Merton problem}
\label{subsec:example 7}
 Let $(X_s)_{t\in [0,T]}$ be the solution of the controlled SDE
$$
dX_s = (\mu \pi_s X_s - c_s)ds + \pi_s \sigma X_s dB_s
$$
started at $X_t = x$, where
$(c_s)_{s\in [0,T]}$
and
$(\sigma_s)_{s\in [0,T]}$ are
square-integrable adapted processes.
We consider the Merton problem
 \begin{equation}
   \label{fjhkdsf}
u(t, x) =
\inf\limits_{(\pi_s)_{t \leq s \leq T}, \ \! ( c_s)_{t \leq s \leq T}}
\mathbb{E}
\left[
\frac{e^{-\rho (T-t)} X_T^{1 - \gamma}}{1 - \gamma}
+
\int_t^T \frac{e^{-\rho (s-t)} c_s^{1 - \gamma}}{1 - \gamma} ds
\right],
\end{equation}
 where $\gamma \in (0,1)$. The solution $u(t,x)$ of \eqref{fjhkdsf}
 satisfies the Hamilton-Jacobi-Bellman (HJB) equation
$$
\partial_t u(t, x) + \sup\limits_{\pi, c}
\left(
(\pi \mu x - c)\partial_x u(t, x) + \frac{\pi^2 \sigma^2 x^2}{2}\partial^2_x u(t, x)
+ \frac{c^{1-\gamma}}{1-\gamma}
\right)
= \rho u(t, x),
$$
 which, by first order condition, can be rewritten as
\begin{equation}
    \label{eq:example 7}
    \partial_t u(t, x)
    - \frac{(\mu \partial_x u(t, x))^2}{2\sigma^2 \partial^2_x u(t, x)}
    + \frac{\gamma}{1 - \gamma} (\partial_x u(t, x))^{1 - 1/\gamma}
    = \rho u(t, x),
\end{equation}
and admits the solution
$$
u(t, x) =
\frac{x^{1-\gamma} (1 + (\alpha - 1) e^{-\alpha (T-t)})^\gamma}
        {\alpha^\gamma (1 - \gamma)},
\qquad (t,x) \in [0,T]\times \real,
$$
where $\alpha := ( 2\sigma^2 \gamma \rho - (1-\gamma)\mu^2 ) / ( 2\sigma^2 \gamma^2 )$.
As the loss function used in the deep Galerkin method uses a
 division by the second derivatives of the neural network function,
see \eqref{eq:dgm loss function} and \eqref{eq:example 7},
it explodes when the second derivatives of
the learned neural network function becomes small
during the training.
 {Hence, in Table~\ref{table:example 7},
 we only present the outputs of the deep branching method
 and of the deep BSDE method of 
 \cite{beck} which deals with second order gradient nonlinearities.  
 Table~\ref{table:example 7} 
 summarizes the results of $10$ independent runs,
 with $\mu = 0.03$, $\sigma = 0.1$, $\gamma = 0.5$, $\rho = 0.01$,
 $T = 0.1$ on the interval $[x_{\rm min},x_{\rm max}]=[100,200]$,
 where we take $M = 10,000$ in the deep branching method.} 
  
\begin{table}[H]
    \centering
    \resizebox{\textwidth}{!}{\begin{tabular}{|l|c||c|c||c|c||c|}
        \hline
        \multicolumn{1}{|c|}{Method} & $d$ & Mean $L^1$-error & Stdev $L^1$-error &
        Mean $L^2$-error & Stdev $L^2$-error & Mean Runtime \\
        \hline
        Deep branching                  & 1 & 8.49E-03 & 7.44E-04 & 1.30E-04 & 2.52E-05 & 54m \\
        \cline{1-1}\cline{3-7}
        Deep BSDE \cite{han2018solving} & 1 & 1.61E+00 & 1.05E-01 & 2.64E+00 & 3.37E-01 & 184m \\
        \hline
\end{tabular}}
		\caption{Summary of numerical results for \eqref{eq:example 7}.}
        \label{table:example 7}
\end{table}

\vskip-0.4cm

\noindent
 An anomaly was detected on the third run when using $\sigma=\sigma_{\tanh}$,
 and it disappeared after 
 changing the activation function to $\sigma=\sigma_{\rm ReLU}$.

\vskip-0.1cm

\begin{figure}[H]
  \centering
 \begin{subfigure}[b]{0.45\textwidth}
    \includegraphics[width=\linewidth]{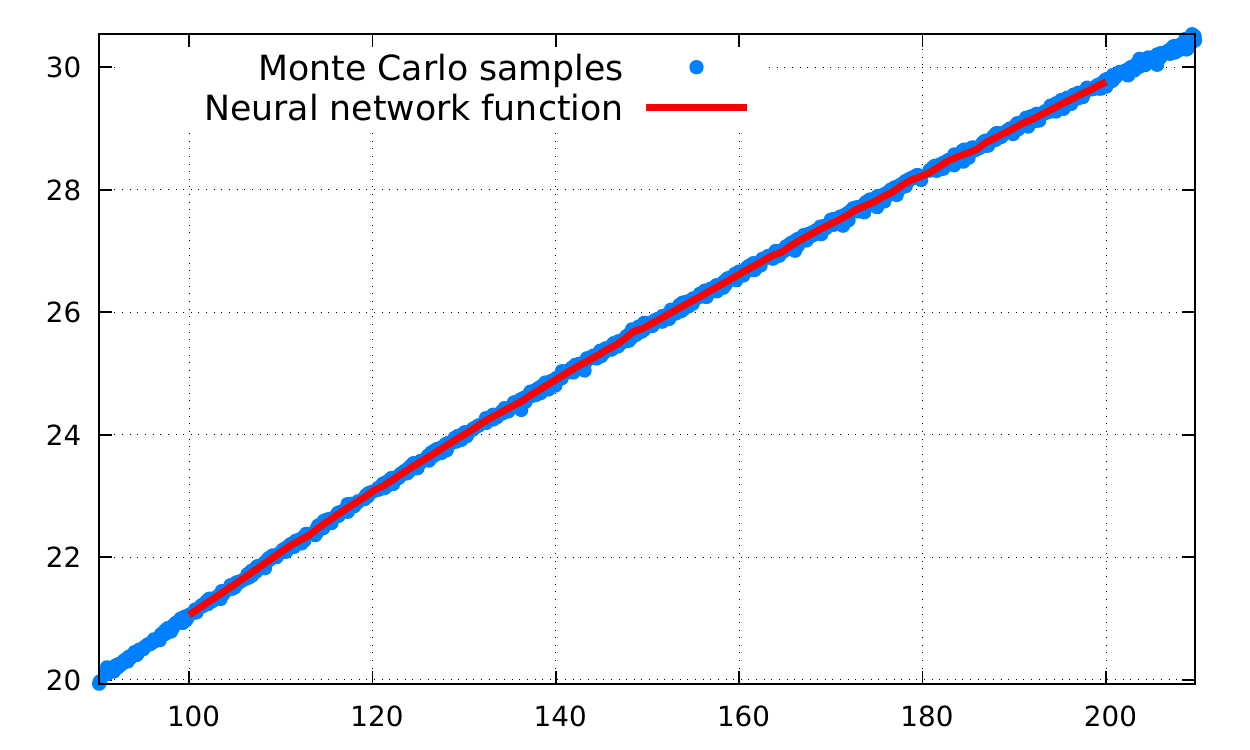}
    \caption{$\sigma=\sigma_{\rm ReLU}$}
 \end{subfigure}
  \begin{subfigure}[b]{0.45\textwidth}
    \includegraphics[width=\linewidth]{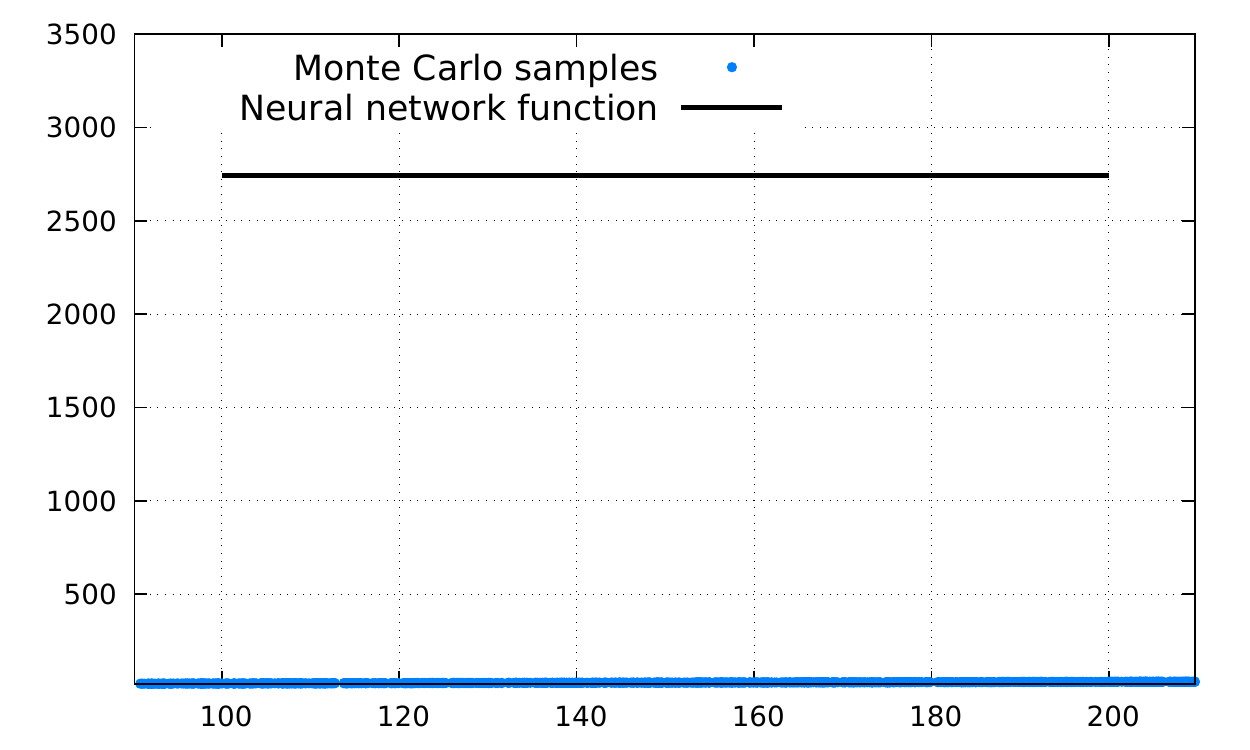}
    \caption{$\sigma=\sigma_{\tanh}$}
 \end{subfigure}
   \vskip-0.2cm
  \caption{Monte Carlo samples and the neural network function on the third run.}
\label{fig:example 7 debug}
\end{figure}

\noindent
In Figure~\ref{fig:example 7 debug}, we plot
the Monte Carlo samples generated by Algorithm~\ref{alg:coding tree}
and the learned neural network function $v(\cdot, \cdot; \theta)$,
see Algorithm~\ref{alg:deep},
for $\sigma=\sigma_{\rm ReLU}$ and for $\sigma=\sigma_{\tanh}$ 
on the third run. 

\vskip-0.3cm

\noindent
Figure~\ref{fig:example 7 debug} shows the consistency, or lack thereof,
between the Monte Carlo samples and the learned neural network function,
which cannot be observed when using the deep Galerkin or deep BSDE method.

  \subsubsection*{e) Third order gradient log nonlinearity}
\label{subsec:example 6}
\begin{figure}[H]
\centering
\includegraphics[width=0.55\textwidth]{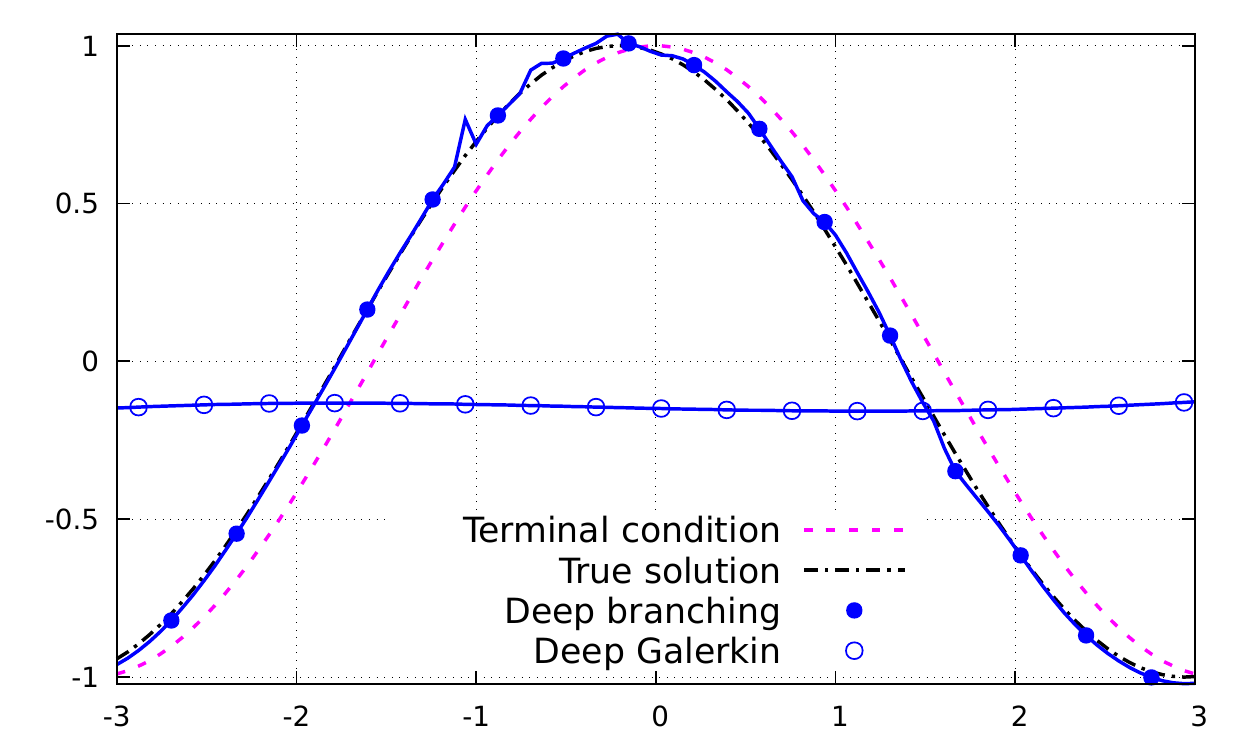}
\vskip-0.1cm
\caption{Deep branching {\em vs} deep Galerkin method 
  for  \eqref{eq:example 6} with $d = 5$ and $T=0.02$.}
\label{fig:example 6}
\end{figure}

{This Example~$e)$ and the next Example~$f)$
 use nonlinearities in terms of third and fourth order gradients, 
 to which the deep BSDE method does not apply. 
 For this reason, comparisons are done only 
 with respect to the Galerkin method.
  }
Consider the equation
\label{p-15}
\begin{equation}
   \label{eq:example 6}
   \partial_t u(t,x) + \frac{\alpha}{d} \sum\limits_{i=1}^d \partial_{x_i} u(t,x)
   + \log \left(
       \frac{1}{d} \sum\limits_{i=1}^d
       \left(\partial^2_{x_i} u(t,x)\right)^2 + \left(\partial^3_{x_i} u(t,x)\right)^2
   \right)
   = 0,
\end{equation}
which admits the solution
$$
u(t,x) = \cos \left( \sum\limits_{i=1}^d x_i + \alpha (T-t)\right),
\qquad (t,x) \in [0,T]\times \real^d.
$$
Table~\ref{table:example 6} summarizes the results of $10$ independent runs,
with $M = 6,000$ in dimension $d = 1$ (resp. $M = 200$ in dimension $d = 5$), $\alpha = 10$,
$T = 0.02$, $x_{\rm min} = -3$, and $x_{\rm max} = 3$.
\begin{table}[H]
    \centering
    \resizebox{\textwidth}{!}{\begin{tabular}{|l|c||c|c||c|c||c|}
        \hline
        \multicolumn{1}{|c|}{Method} & $d$ & Mean $L^1$-error & Stdev $L^1$-error &
        Mean $L^2$-error & Stdev $L^2$-error & Mean Runtime \\
        \hline
        Deep branching                        & 1 & 5.82E-03 & 1.27E-03 & 5.52E-05 & 2.01E-05 & 78m \\
        \cline{1-1}\cline{3-7}
        Deep Galerkin \cite{sirignano2018dgm} & 1 & 7.50E-02 & 3.15E-02 & 9.00E-03 & 7.34E-03 & 83m \\
                \hline
        Deep branching                        & 5 & 2.77E-02 & 1.13E-02 & 3.52E-03 & 4.65E-03 & 183m \\
        \cline{1-1}\cline{3-7}
Deep Galerkin \cite{sirignano2018dgm} & 5 & 6.38E-01 & 5.74E-03 & 5.18E-01 & 1.08E-02 & 369 \\
                \hline
\end{tabular}}
		\caption{Summary of numerical results for \eqref{eq:example 6}.}
        \label{table:example 6}
\end{table}
\vskip-0.3cm
\noindent
In the case of log nonlinearity with a third order gradient term,
our method appears more
accurate than the deep Galerkin method in dimensions $d=1$ and $d=5$.
 Figure~\ref{fig:example 6} presents a numerical comparison
 on the average performance of 10 runs.

\vskip-0.1cm

\vskip-0.3cm

\subsubsection*{f) Fourth order gradient cosine nonlinearity}
\label{subsec:example 5}
Consider the equation
 \begin{equation}
   \label{eq:example 5}
   \partial_t u(t,x) + \frac{\alpha}{d} \sum\limits_{i=1}^d \partial_{x_i} u(t,x)
   + u(t,x) - \left(\frac{\Delta u(t,x)}{12d} \right)^2
   + \frac{1}{d} \sum\limits_{i=1}^d
        \cos \left( \frac{\pi \partial^4_{x_i} u(t,x)}{4!} \right)
   = 0,
\end{equation}
 which admits the solution
 $$
 u(t,x)=\varphi \left( \sum\limits_{i=1}^d x_i + \alpha ( T-t)\right),
 \qquad (t,x) \in [0,T]\times \real^d,
 $$
 where
 $\varphi (y) := y^4 + y^3 + by^2 + cy + d$ for $y \in \real$,
 $b = -36/47$, $c = 24b$, $d = 4b^2$, and $\alpha = 10$.

\medskip
Table~\ref{table:example 5} summarizes the results of $10$ independent runs,
with $M = 2,500$ in dimension $d = 1$ (resp. $M = 50$ in dimension $d = 5$), $\alpha = 10$, $T = 0.04$, $x_{\rm min} = -5$, and $x_{\rm max} = 5$.

\begin{figure}[H]
\centering
\includegraphics[width=0.55\textwidth]{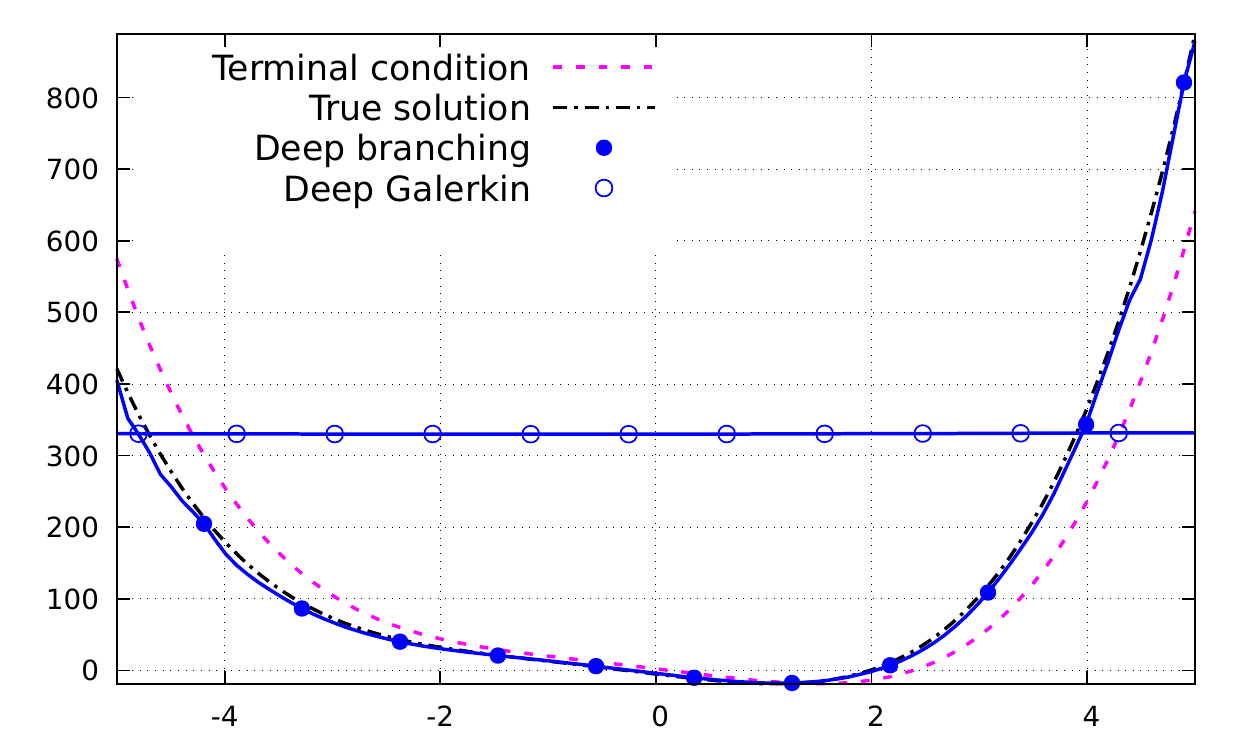}
\vskip-0.1cm
\caption{Deep branching {\em vs} deep Galerkin method for
  \eqref{eq:example 5} with $d = 5$ and $T=0.04$.}
\label{fig:example 5}
\end{figure}

\begin{table}[H]
    \centering
    \resizebox{\textwidth}{!}{\begin{tabular}{|l|c||c|c||c|c||c|}
        \hline
        \multicolumn{1}{|c|}{Method} & $d$ & Mean $L^1$-error & Stdev $L^1$-error &
        Mean $L^2$-error & Stdev $L^2$-error & Mean Runtime \\
        \hline
        Deep branching                        & 1 & 9.62E+00 & 1.50E+00 & 3.76E+02 & 1.62E+02 & 128m \\
        \cline{1-1}\cline{3-7}
        Deep Galerkin \cite{sirignano2018dgm} & 1 & 2.81E+01 & 2.77E+01 & 2.31E+03 & 4.45E+03 & 146m \\
        \hline
        Deep branching                        & 5 & 1.01E+01 & 1.16E+00 & 3.49E+02 & 1.62E+02 & 259m \\
        \cline{1-1}\cline{3-7}
Deep Galerkin \cite{sirignano2018dgm} & 5 & 2.57E+02 & 1.18E+00 & 7.75E+04 & 6.55E+02 & 670 \\
        \hline
\end{tabular}}
		\caption{Summary of numerical results for \eqref{eq:example 5}.}
        \label{table:example 5}
\end{table}
\vskip-0.3cm
\noindent
In the case of cosine nonlinearity with a fourth order gradient,
our method appears more
accurate than the deep Galerkin methods in dimensions $d=1$ and $d=5$.
 Figure~\ref{fig:example 5} presents a numerical comparison
 on the average performance of 10 runs.

\vskip-0.1cm

\appendix

\section{Multidimensional extension}
In this section we sketch the argument extending
Theorem~1 in \cite{penent2022fully} to the multidimensional case,
 and leading to \eqref{eq:feynman kac}.
 For this,
 given $g \in C^{0, \infty}([0, T] \times \real^d)$
 and $\mu \in \inte^n$ such that
 $|\mu| \geq 1$, we will use the multivariate Fa\`a di Bruno formula
\begin{equation} 
  \label{faadibruno}
  \partial_\mu g^*(u)(t, x) =
 \left(\prod_{i = 1}^d \mu_i! \right)
    \sum\limits_{\substack{1 \le \nu_1 + \cdots + \nu_n \le \abs{\mu} \\
                           1 \le s \le \abs{\mu}}}
   (\partial_{\nu} g)^* (u)(t,x)
    \hskip-0.54cm
    \sum\limits_{\substack{1 \le \abs{k_1}, \dots, \abs{k_s}, \
                           0 \prec l^1 \prec \cdots \prec l^s \\
                           k^i_1 + \cdots + k^i_s = \nu_i, \
                           i = 1, \dots, n \\
                           \abs{k_1}l_j^1 + \cdots + \abs{k_s}l_j^s = \mu_j, \
                           j = 1, \dots, d
                           }}
    \prod_{\substack{1 \le r \le s \\
                     1 \le q \le n
                    }}
    \frac{(\partial_{l^r + \lambda^q} u (t,x))^{k_r^q}}{k_r^q!
                \left(l_1^r! \cdots l_d^r!\right)^{k_r^q}},
\end{equation} 
 see Theorem~2.1 in \cite{constantine}, applied to the function
 $g^*(u)(t, x) := 
 g\big(\partial_{\lambda^1}u(t,x),\ldots , \partial_{\lambda^n}u(t,x)\big)$.

\medskip
\noindent
 We have
\begin{eqnarray*}
  \lefteqn{
    \partial_t g^*(u) + \frac{1}{2}\Delta g^*(u)
  }
  \\
  & = &
  \sum\limits_{p = 1}^n
      \partial_{\lambda^p} \left(\partial_t u + \frac{1}{2}\Delta u \right)
      (\partial_{\bm{1}_p} g)^*(u)
  + \frac{1}{2} \sum\limits_{i = 1}^n \sum\limits_{j = 1}^n \sum\limits_{k = 1}^d
        \left(\partial_{\lambda^i + \bm{1}_k} u\right)
        \left(\partial_{\lambda^j + \bm{1}_k} u\right)
        (\partial_{\bm{1}_i + \bm{1}_j} g)^*(u)
  \\
  &= &
  - \sum\limits_{p = 1}^n
      (\partial_{\bm{1}_p} g)^*(u)
  \partial_{\lambda^p} f^*(u)
  + \frac{1}{2} \sum\limits_{i = 1}^n \sum\limits_{j = 1}^n \sum\limits_{k = 1}^d
        \left(\partial_{\lambda^i + \bm{1}_k} u\right)
        \left(\partial_{\lambda^j + \bm{1}_k} u\right)
        (\partial_{\bm{1}_i + \bm{1}_j} g)^*(u)
  \\
  &= &
    - \sum\limits_{p = 1}^n
  \bm{1}_{\{ \lambda^p = 0\}} (\partial_{\bm{1}_p} g)^*(u)
   f^*(u)
  \\
   & & 
  - \sum\limits_{p = 1}^n
      (\partial_{\bm{1}_p} g)^*(u)
    \left(\prod_{i = 1}^d \lambda^p_i! \right)
    \sum\limits_{\substack{1 \le \nu_1 + \cdots + \nu_n \le \abs{\lambda^p} \\
                           1 \le s \le \abs{\lambda^p}}}
   (\partial_{\nu} f)^*(u)
    \sum\limits_{\substack{1 \le \abs{k_1}, \dots, \abs{k_s}, \
                           0 \prec l^1 \prec \cdots \prec l^s \\
                           k^i_1 + \cdots + k^i_s = \nu_i, \
                           i = 1, \dots, n \\
                           \abs{k_1}l_j^1 + \cdots + \abs{k_s}l_j^s = \lambda^p_j, \
                           j = 1, \dots, d
                           }}
    \prod_{\substack{1 \le r \le s \\
                     1 \le q \le n
                    }}
    \frac{(\partial_{l^r + \lambda^q} u)^{k_r^q}}{k_r^q!
                \left(l_1^r! \cdots l_d^r!\right)^{k_r^q}}
  \\
  & & \qquad \qquad \qquad
      + \frac{1}{2} \sum\limits_{i = 1}^n \sum\limits_{j = 1}^n \sum\limits_{k = 1}^d
        \left(\partial_{\lambda^i + \bm{1}_k} u\right)
        \left(\partial_{\lambda^j + \bm{1}_k} u\right)
        (\partial_{\bm{1}_i + \bm{1}_j} g)^*(u).
\end{eqnarray*}
 Rewriting the above equation in integral form yields
 \begin{align}
   \nonumber
   & g^*(u)(t,x) = \int_{\real^d} \varphi (T-t,y-x) g ( \phi(y) ) dy
   \\
   \nonumber
   &
   + \int_t^T \int_{\real^d} \varphi (s-t,y-x)
    \\
  \nonumber
    &
   \biggl(
  \sum\limits_{p = 1}^n
  \bm{1}_{\{ \lambda^p = 0\}} (\partial_{\bm{1}_p} g)^*(u)
   f^*(u) 
        -\frac{1}{2} \sum\limits_{i = 1}^n \sum\limits_{j = 1}^n \sum\limits_{k = 1}^d
        \left(\partial_{\lambda^i + \bm{1}_k} u(s,y)\right)
        \left(\partial_{\lambda^j + \bm{1}_k} u(s,y)\right)
        (\partial_{\bm{1}_i + \bm{1}_j} g)^*(u)
    \\
  \nonumber
    & 
   + \sum\limits_{p = 1}^n
       (\partial_{\bm{1}_p} g)^*(u)
    \left(\prod_{i = 1}^d \lambda^p_i! \right)
    \sum\limits_{\substack{1 \le \nu_1 + \cdots + \nu_n \le \abs{\lambda^p} \\
                           1 \le s \le \abs{\lambda^p}}}
   (\partial_{\nu} f)^* (u) 
    \sum\limits_{\substack{1 \le \abs{k_1}, \dots, \abs{k_s}, \
                           0 \prec l^1 \prec \cdots \prec l^s \\
                           k^i_1 + \cdots + k^i_s = \nu_i, \
                           i = 1, \dots, n \\
                           \abs{k_1}l_j^1 + \cdots + \abs{k_s}l_j^s = \lambda^p_j, \
                           j = 1, \dots, d
                           }}
    \prod_{\substack{1 \le r \le s \\
                     1 \le q \le n
                    }}
    \frac{(\partial_{l^r + \lambda^q} u(s,y))^{k_r^q}}{k_r^q!
                \left(l_1^r! \cdots l_d^r!\right)^{k_r^q}}
   \biggr) dy ds. 
   \\
   \label{refg1}
 \end{align}
 Similarly, for $\mu \in \mathbb{N}^d$,
 by the Fa\`a di Bruno formula \eqref{faadibruno} 
 we have 
\begin{align}
  \nonumber
  &
 \partial_\mu u(t,x) =
 \int_{\real^d} \varphi (T-t,y-x)
  \partial_\mu u (T,y)
  dy
 \\
\nonumber
&
\hskip-0.2cm
+
\int_t^T \int_{\real^d} 
 \sum_{
     \footnotesize \substack{
               1 \le \nu_1 + \cdots + \nu_n \le \abs{\mu} 
               \\
               1 \le s \le \abs{\mu}
 }}
 \sum_{
     \footnotesize \substack{
               1 \le \abs{k_1}, \dots, \abs{k_s}, \
               0 \prec l^1 \prec \cdots \prec l^s \\
               k^i_1 + \cdots + k^i_s = \nu_i, \
               i = 1, \dots, n \\
               \abs{k_1}l_j^1 + \cdots + \abs{k_s}l_j^s = \mu_j, \
               j = 1, \dots, d
                }}
\frac{\prod\limits_{i = 1}^d \mu_i! }
     {\prod\limits_{\footnotesize \substack{1 \le r \le s \\
                     1 \le q \le n}}
      k_r^q! \left(l_1^r! \cdots l_d^r!\right)^{k_r^q}
     }
     (\partial_{\nu} f)^*(u)
     \sum_{\footnotesize \substack{1 \le r \le s \\
                     1 \le q \le n
                    }}
     \big(
     \partial_{l^r + \lambda^q} u (s,y) 
     \big)^{k_r^q} dy ds. 
   \\
   \label{refg2}
\end{align}
Combining \eqref{refg1} and \eqref{refg2} yields the equation
\begin{equation}
 \label{s1}
 c(u)(t,x) = \int_{-\infty}^\infty \varphi (T-t,y-x) c ( u)(T,y) dy
+
\sum_{Z \in \mathcal{M}(c)}
\int_t^T \int_{-\infty}^\infty \varphi (s-t,y-x)
\prod_{z \in Z}  z(u)(s,y) dy ds,
\end{equation}
$(t,x)\in [0,T]\times \real$,
 for any code $c\in \mathcal{C}$, as in Lemma~2.3 of \cite{penent2022fully}.
The dimension-free argument of Theorem~1 in \cite{penent2022fully}
then shows that \eqref{eq:feynman kac} holds provided
that ${\cal H} ({t, x, \rm Id} )$ is integrable and the solution of
\eqref{s1} is unique.

\footnotesize

\newcommand{\etalchar}[1]{$^{#1}$}
\def\cprime{$'$} \def\polhk#1{\setbox0=\hbox{#1}{\ooalign{\hidewidth
  \lower1.5ex\hbox{`}\hidewidth\crcr\unhbox0}}}
  \def\polhk#1{\setbox0=\hbox{#1}{\ooalign{\hidewidth
  \lower1.5ex\hbox{`}\hidewidth\crcr\unhbox0}}} \def\cprime{$'$}

\end{document}